\documentclass[12pt]{scrartcl}

\usepackage{amsmath,amssymb,mathtools}

\usepackage{upgreek}

\usepackage{lineno}

\newcommand{\bfsf}[1]{\boldsymbol{\mathsf{#1}}}

\newcounter{thm}
\newenvironment{thm}[1][]{\refstepcounter{thm}\par\medskip
   \noindent \textbf{\sffamily Theorem~\thethm. #1} \slshape }{\medskip}
   
   \newcounter{propo}

\newcounter{rem}
\newenvironment{rem}[1][]{\refstepcounter{rem}\par\medskip
   \noindent \textbf{\sffamily Remark~\therem. #1} \rmfamily }{\medskip}
   
   \newcounter{defi}

\mathtoolsset{showonlyrefs=true}

\makeatletter
\g@addto@macro\bfseries{\boldmath}
\makeatother

\usepackage{amsmath,amssymb,amsthm,mathtools,array}
\usepackage{graphicx}

\usepackage{lineno}

\mathtoolsset{showonlyrefs=true}

 \newtheorem{cor}[thm]{Corollary}

 \theoremstyle{definition}
 
 \theoremstyle{remark}

 \numberwithin{equation}{section}

\usepackage{tikz}

\usetikzlibrary{patterns}

\newlength{\hatchspread}
\newlength{\hatchthickness}
\newlength{\hatchshift}
\newcommand{\hatchcolor}{}
\tikzset{hatchspread/.code={\setlength{\hatchspread}{#1}},
         hatchthickness/.code={\setlength{\hatchthickness}{#1}},
         hatchshift/.code={\setlength{\hatchshift}{#1}},
         hatchcolor/.code={\renewcommand{\hatchcolor}{#1}}}
\tikzset{hatchspread=3pt,
         hatchthickness=0.4pt,
         hatchshift=0pt,
         hatchcolor=black}
\pgfdeclarepatternformonly[\hatchspread,\hatchthickness,\hatchshift,\hatchcolor]
   {custom north west lines}
   {\pgfqpoint{\dimexpr-2\hatchthickness}{\dimexpr-2\hatchthickness}}
   {\pgfqpoint{\dimexpr\hatchspread+2\hatchthickness}{\dimexpr\hatchspread+2\hatchthickness}}
   {\pgfqpoint{\dimexpr\hatchspread}{\dimexpr\hatchspread}}
   {
    \pgfsetlinewidth{\hatchthickness}
    \pgfpathmoveto{\pgfqpoint{0pt}{\dimexpr\hatchspread+\hatchshift}}
    \pgfpathlineto{\pgfqpoint{\dimexpr\hatchspread+0.15pt+\hatchshift}{-0.15pt}}
    \ifdim \hatchshift > 0pt
      \pgfpathmoveto{\pgfqpoint{0pt}{\hatchshift}}
      \pgfpathlineto{\pgfqpoint{\dimexpr0.15pt+\hatchshift}{-0.15pt}}
    \fi
    \pgfsetstrokecolor{\hatchcolor}
    \pgfusepath{stroke}
   }

\mathtoolsset{showonlyrefs=true}

\makeatletter
\g@addto@macro\bfseries{\boldmath}
\makeatother

\begin{document}

\begin{center}
\bfseries{\sffamily\Large On the Clifford Algebraic Description of the\\ \smallskip
Geometry of a 3D Euclidean Space}
\end{center}

\bigskip
\begin{center}
{\bfseries\large Jayme Vaz Jr.}$^{a,}$\footnote{\texttt{<vaz@ime.unicamp.br>}}
{\bfseries\large and} {\bfseries\large Stephen Mann}$^{b,}$\footnote{\texttt{<smann@uwaterloo.ca>}}\\
\medbreak
\small 
$^a$ Department of Applied Mathematics\\
University of Campinas\\ 
Campinas, SP, Brazil\\
\smallskip
$^b$ Cheriton School of Computer Science\\
University of Waterloo\\
Waterloo, ON, Canada 
\end{center}
\normalsize

\bigskip

\footnotesize 

\begin{center}
\begin{minipage}{11 cm}
\textbf{Abstract:} 
We discuss how transformations in a three dimensional euclidean space
can be described in terms of the Clifford algebra $\mathcal{C}\ell_{3,3}$
of the quadratic space $\mathbb{R}^{3,3}$. 
We show that this algebra describes in a 
unified way the operations of reflection, rotations (circular 
and hyperbolic), translation, shear and non-uniform 
scale. Moreover, using the concept of Hodge duality, 
we define an operation called cotranslation, and show
that the operation of perspective projection can be 
written in this Clifford algebra as a composition of the 
translation and cotranslation
operations.  We also show that the operation of pseudo-perspective can be
implemented using the cotranslation operation. An important point 
is that the expression for the operations of reflection and rotation in 
$\mathcal{C}\ell_{3,3}$ preserve the subspaces that
can be associated with the algebras 
$\mathcal{C}\ell_{3,0}$ and $\mathcal{C}\ell_{0,3}$, so that
reflection and rotation can be expressed in terms of $\mathcal{C}\ell_{3,0}$ or 
$\mathcal{C}\ell_{0,3}$, as well-known. However, all other operations 
mix those subspaces in such a way that they need to be
expressed in terms of 
the full Clifford algebra $\mathcal{C}\ell_{3,3}$.  
An essential aspect of our formulation is the representation 
of points in terms of objects called paravectors. Paravectors
have been used previously to represents points 
in terms of an algebra closely related
to the Clifford algebra $\mathcal{C}\ell_{3,3}$. We compare 
these different approaches.
\end{minipage}
\end{center}

\normalsize

\bigskip

\section{Introduction}

There is a deep connection between geometry and algebra, 
and exploiting this connection usually benefits 
the studies and advances of them. Computer graphics is one of many 
areas where geometry and algebra play 
a key role. The concepts
of affine spaces and projective spaces are fundamental 
for the theoretical basis of computer graphics, 
and linear algebra is a fundamental tool in computer graphics
for geometric computations. Nevertheless, the success of this 
formalism should not be seen as a hindrance to the study of 
the use of other algebraic systems. In fact, quaternions~\cite{quaternions,Goldman}
were shown to be an efficient tool for the
interpolation of sequences of orientations in computer graphics, 
and other studies~\cite{DFM} involving applications of other algebraic systems 
in computer graphics have appeared in the
literature, together with studies
trying to identify geometric spaces that may be
more suitable as an ambient space for 
computer graphics~\cite{Goldman2000,Goldman2002}. 
Models based on Clifford algebras~\cite{GS,Dorst2,DGM17} have
been proposed as an alternative algebraic framework 
for computer graphics. 
 
Recently we have proposed a new model for the
description of a geometrical space 
based on the exterior algebra of a vector space~\cite{VM2018}. 
In this model, points are described by objects called
paravectors, line segments are described by biparavectors, 
plane fragments are described by triparavectors, and so on 
for higher dimensional spaces. A $k$-{\em paravector} is a
sum of a $(k-1)$-vector and a $k$-vector, which are
elements of the exterior algebra. In~\cite{VM2018} we
exploited the algebra of $k$-paravectors 
to describe some properties of a geometrical space
relevant for computer graphics, in particular we 
studied geometric transformations of points, lines and
planes. These transformations were studied in 
terms of an algebra of transformations that is analogous to the 
algebra of creation and annihilation operators in quantum 
theory. We showed that this model describes in a 
unified way the operations of reflection, rotations (circular 
and hyperbolic), translation, shear and non-uniform 
scale; and, using the concept of Hodge duality, 
we defined an operation called {\em cotranslation}, and showed
that the operation of perspective projection can be 
written as a composition of the translation and cotranslation
operations, and that the operation of pseudo-perspective can be
implemented using the cotranslation operation.
The advantages of our model are that it contains points,
line segments, and plane sectors in a natural way, and it includes all
affine and projective transformations.
In relation to these transformations, the derivations of shear, 
non-uniform scaling, perspective and
pseudo-perspective are easy; on the other hand, the disadvantages  
include that the derivation of reflection is more complicated
than in other approaches, as well as the derivations of rotation and
hyperbolic rotation.

There is a deep relationship between
the algebra of creation and annihilation operators of fermions 
in quantum mechanics and the Clifford algebra 
of a quadratic space~\cite{VazRocha}. Since our exterior algebra
model~\cite{VM2018} is based on an algebra analogous
to the algebra of creation and annihilation operators, 
this relationship leads us to consider the formulation of 
our model in terms of a Clifford algebra. 
In addition to the interest in the problem itself, in 
particular for the comparison of the usefulness of 
the different algebraic frameworks, 
this Clifford algebra formulation would also allow us to be able to compare in a 
more direct manner our model with other models 
in the literature based on Clifford algebras. 
These are the main purposes of this work. 

We have organized this paper as follows. 
In section~\ref{sec.2} we briefly discuss the 
Clifford algebra of the Euclidean space $V_3 \simeq \mathbb{R}^3$. 
We can associate to this space two different 
Clifford algebras, usually denoted by $\mathcal{C}\ell_{3,0}$ 
and $\mathcal{C}\ell_{0,3}$.
Although these two algebras are not isomorphic, 
they are
equivalent in the sense that geometric 
transformations can be equivalently described
by these two algebras. However, the set 
of geometric transformations described by
these two algebras is limited and has to 
be enlarged if we want an algebraic 
framework that can handle transformations
used in computer graphics. In section~\ref{sec.3} 
we discuss the Clifford algebra $\mathcal{C}\ell_{3,3}$ 
and interpret this algebra as being associated with
the vector space $V_3$ and its dual $V_3^\ast$. 
The key point here is how we identify the vector
space $V_3$ inside the algebra $\mathcal{C}\ell_{3,3}$, 
and this identification underpins our model. 
Then in section~\ref{vect_point} we briefly 
discuss the ideas of~\cite{VM2018} concerning
the use of paravectors to represent 
points, and re-express these ideas in terms of $\mathcal{C}\ell_{3,3}$. 
In section~\ref{Trans.points} we exploit this
representation to study some transformations
of points. We will see that reflections and
rotations leave invariant the subspaces that
can be associated with the algebras 
$\mathcal{C}\ell_{3,0}$ and $\mathcal{C}\ell_{0,3}$, 
but other transformations like shear or translation 
do not leave these subspaces invariant. This is
why, from the point of view of our model, 
reflection and rotation can be (nicely) 
described by $\mathcal{C}\ell_{3,0}$ and $\mathcal{C}\ell_{0,3}$, 
but not transformations like shear or non-uniform 
scale transformation. 
We also define a transformation called 
cotranslation and show how cotranslation can be used
to describe the operations of perspective projection and 
pseudo-perspective projection. 
In section~\ref{sec.gen.trans} we study 
the general form of a transformation in 
our algebraic framework and conclude that 
it is a projective transformation. 
 In section~\ref{sec.7}
we present our conclusions.

\section{The 3D Euclidean space and the Clifford algebras $\mathcal{C}\ell_{3,0}$ and $\mathcal{C}\ell_{0,3}$}
\label{sec.2}

In what follows we restrict
things to three dimensions and to the Euclidean case 
because this case is of special interest to computer graphics.
Consider a tridimensional vector space $V_3$ and its
dual vector space $V_3^\ast$. 
Although $V_3 \simeq V_3^\ast \simeq \mathbb{R}^3$,  
we must resist the temptation to identify these spaces. 
This separation is clear from the mathematical point of view, 
but often, in the applications of these concepts 
to particular problems, these two spaces 
are shuffled and the distinction ends up being lost. 
The notation we use is one
tool to manifest the distinction between 
$V_3$ and $V_3^\ast$. 
Let $\mathfrak{B} = \{\vec{e}_1,\vec{e}_2,\vec{e}_3\}$ 
be a basis for $V_3$ and $\mathfrak{B}^\ast = \{\theta^1,\theta^2,\theta^3\}$ 
be a basis for $V_3^\ast$. The basis $\mathfrak{B}^\ast$ is the
dual basis of $\mathfrak{B}$ when $\theta^i(\vec{e}_j) = \delta^i_j$. 
Then, for a vector $\vec{v} = v^i \vec{e}_i$ (summation convention 
is implied) and a covector $\alpha = \alpha_i \theta^i$, we have 
$\alpha(\vec{v}) = \alpha_i v^i$. 
Given a metric $g: V \times V \rightarrow \mathbb{R}$, 
we  write $g(\vec{v},\vec{u}) = g_{ij} v^i u^j$, where 
$g_{ij} = g(\vec{e}_i,\vec{e}_j) = g_{ji}$ 
with $\vec{v} = v^i\vec{e}_i$ and $\vec{u} = u^i\vec{e}_i$.  
An object like $g_{ij}\theta^j$ is a linear combination 
of covectors, and as such $g_{ij}\theta^j \in V_3^\ast$. 
If we want to keep clear the
distinction between vectors and covectors, the notation has to 
reflect this, and we will denote $g_{ij}\theta^j = \vec{e}_i{}^\ast$. 
In other words, $\vec{e}_i$ is a vector, while 
$\vec{e}_i{}^\ast$ is a covector. 

Consider the problem of multiplying vectors. 
The real Clifford algebra $\mathcal{C}\ell_{3,0}$ of the
Euclidean 3-space $(V_3,g)$ is the associative
algebra generated by $\{ \gamma(\vec{v}) \vert 
\vec{v} \in V_3\}$
and $\{a 1 \vert a \in \mathbb{R}\}$, 
where $\gamma: V_3 \rightarrow \mathcal{C}\ell_{3,0}$, 
subject to the
condition~\cite{DFM,VazRocha,Lounesto}
\begin{equation}
  \gamma(\vec{v}) \gamma(\vec{u})+ \gamma(\vec{u}) 
  \gamma(\vec{v}) = 2 g ( \vec{v} , \vec{u} ) , \label{eq:pro_cli1}
\end{equation}
where the Clifford product is denoted by juxtaposition. 

A completely equivalent definition of the 
Clifford algebra of the Euclidean 3-space $(V_3,g)$, 
which we will denote for the moment by $\mathcal{C}\ell^\prime_{3,0}$, 
is that of 
the associative
algebra generated by $\{ \gamma(\vec{v}) \vert \vec{v} \in V_3\}$
and $\{a 1 \vert a \in \mathbb{R}\}$, where $\gamma: V_3 \rightarrow 
\mathcal{C}\ell_{3,0}^\prime$, subject to the
condition~\cite{DFM,VazRocha,Lounesto}
\begin{equation}
  \gamma(\vec{v}) \gamma(\vec{u}) + \gamma(\vec{u})
  \gamma(\vec{v}) =  - 2 g ( \vec{v} , \vec{u} ) . \label{eq:pro_cli1.alt}
\end{equation}
The difference between this definition and the former one is the
 sign on the right hand side of eq.\eqref{eq:pro_cli1.alt} in 
relation to the right hand side of eq.\eqref{eq:pro_cli1}. This second definition is the 
definition adopted by several authors as~\cite{Porteous,Harvey}. 

It is clear that this latter definition of the Clifford algebra of the space 
$(V_3,g)$ is the same of the first definition for the case of 
the Clifford algebra of the space $(V_3,g^\prime = -g)$. So we will 
abandon the notation $\mathcal{C}\ell_{3,0}^\prime$ 
and the Clifford algebra subject to eq.\eqref{eq:pro_cli1.alt} by $\mathcal{C}\ell_{0,3}$.

Therefore, given a vector space $V_3$, we can 
associate two different Clifford algebras to it, one
using the Clifford map as in eq.\eqref{eq:pro_cli1}, 
which we will denote by 
$\gamma^+: V_3 \rightarrow \mathcal{C}\ell_{3,0}$, 
and another using the Clifford map  as in eq.\eqref{eq:pro_cli1.alt}, 
which we will denote by 
$\gamma^-: V_3 \rightarrow \mathcal{C}\ell_{0,3}$.

The algebras $\mathcal{C}\ell_{3,0}$ and $\mathcal{C}\ell_{0,3}$ 
are not isomorphic; indeed, from the representation theory~\cite{VazRocha} 
of Clifford algebras, we know that $\mathcal{C}\ell_{3,0} \simeq 
\mathcal{M}(2,\mathbb{C})$, where $\mathcal{M}(2,\mathbb{C})$ 
is the algebra of $2\times 2$ complex matrices, and 
$\mathcal{C}\ell_{0,3} \simeq \mathbb{H}\oplus \mathbb{H}$, 
where $\mathbb{H}$ denotes the algebra of quaternions. 
However, although $\mathcal{C}\ell_{3,0}$ and 
$\mathcal{C}\ell_{0,3}$ are not \textit{algebraically} 
equivalent, they are \textit{geometrically} equivalent~\cite{Botman}. 

Let us formulate the $\mathcal{C}\ell_{3,0}$ and $\mathcal{C}\ell_{0,3}$
algebras in some details. 
We denote 
\begin{equation}
\mathbf{e}_i = \gamma(\vec{e}_i) , \qquad i = 1,2,3. 
\end{equation}
Then, from eq.\eqref{eq:pro_cli1} and eq.\eqref{eq:pro_cli1.alt}, we have
\begin{equation}
(\mathbf{e}_i)^2 = \eta_i = \pm 1, \qquad i = 1,2,3,  
\end{equation}
where $\eta_i = 1$ in the case of eq.\eqref{eq:pro_cli1} or 
$\eta = -1$ in the case of eq.\eqref{eq:pro_cli1.alt}. 
 Then an arbitrary element of $\mathcal{C}\ell_{3,0}$ 
or $\mathcal{C}\ell_{0,3}$ 
is  a sum of elements of the form
\begin{equation}
\label{grad}
(\mathbf{e}_1)^{\mu_1}(\mathbf{e}_2)^{\mu_2}
(\mathbf{e}_3)^{\mu_3} , 
\end{equation}
with $\mu_i = 0$ or $1$ ($i = 1,2,3$). 
The unit of the algebra corresponds to $\mu_1 = \mu_2 = \mu_3 = 0$. 
There is a $Z_n$-gradation
and we call an element in eq.\eqref{grad} a $k$-vector
according to $\mu_1 + \mu_2 + \mu_3 = k$. A $0$-vector is 
a scalar; a $1$-vector is a vector associated
with a class of equipollent line segments; a $2$-vector is 
a vector associated with the equivalence class of oriented 
plane fragments with the same area and direction; and a 
$3$-vector is a vector associated with the equivalence 
class of volume elements with the same orientation and 
volume. The 3-vector
$\mathbf{e}_{123} = \mathbf{e}_1\mathbf{e}_2\mathbf{e}_3$
has a special role since it belongs to the center
of $\mathcal{C}\ell_{3,0}$ or $\mathcal{C}\ell_{0,3}$. 
The sum of a scalar and a vector is a {\em paravector}, 
representing a weighted point \cite{VM2018}. Then the paravector 
representing an affine point with positive orientation 
is an element of the form $\boldsymbol{\mathsf{P}} = 1 + \mathbf{p}$, 
where $\mathbf{p} = p^i \mathbf{e}_i$. 

We denote the vector space of $k$-vectors by
$\bigwedge_k(V_3) = \bigwedge_k$
and assume the convention that $\bigwedge_0 = \mathbb{R}$
and $\bigwedge_1 = V_3$.
The projectors $\bigwedge = \oplus_{j=0}^3 \bigwedge_j
\rightarrow \bigwedge_k$ are denoted by
$\langle \hspace{1ex}\rangle_k$. The 
space of paravectors will be denoted by 
$\prod(V_3) = \bigwedge^0 \oplus \bigwedge^1$. 
An arbitrary element of $\bigwedge = \oplus_{j=0}^3 \bigwedge_j $ 
is called a \textit{multivector}. 
We also employ the notation $\mathbf{e}_{ij} =
\mathbf{e}_i \mathbf{e}_j$ for $i \neq j$.

The Clifford product of a vector and an arbitrary
$k$-vector $A_k$ can be written, independently of the
signature of the underlying space, as
\begin{equation}
\label{clifford product}
\mathbf{v} A_k =
\mathbf{v} \cdot A_k + \mathbf{v} \wedge A_k  ,
\end{equation}
where the interior and exterior products are defined, respectively, as
\begin{equation}
\begin{split}
  &\mathbf{v} \cdot A_k = \langle \mathbf{v} A_k \rangle_{k-1} = \frac{1}{2} \left( \mathbf{v} A_k
    - (-1)^k A_k \mathbf{v} \right)  , \\
  &\mathbf{v} \wedge A_k = \langle \mathbf{v} A_k \rangle_{k+1} = \frac{1}{2} \left(
    \mathbf{v} A_k + (-1)^k A_k \mathbf{v} \right)  .
\end{split}
\end{equation}
An important relation between these two products is 
\begin{equation}
\label{Leibniz}
\mathbf{v}\cdot (A_k \wedge B_j) = 
(\mathbf{v}\cdot A_k)\wedge B_j + (-1)^k A_k\wedge (\mathbf{v}\cdot B_j) . 
\end{equation}

We can also define three operations in $\mathcal{C}\ell_{3,0}$,  
or $\mathcal{C}\ell_{0,3}$, 
called graded involution (or parity, denoted by a hat), reversion
(denoted by a tilde) and (Clifford) conjugation (denoted
by a bar), respectively as
\begin{equation}
\label{operations}
\hat{A}_k = (-1)^k A_k , \qquad
\tilde{A}_k = (-1)^{k(k-1)/2}A_k , \qquad
\bar{A}_k = \tilde{\hat{A}}_k = \hat{\tilde{A}}_k .
\end{equation}
The reversion operation receives its name from the property that
\begin{equation}
\widetilde{AB} = \tilde{B}\tilde{A} .
\end{equation}
Conjugation is the composition of graded involution 
and reversion. From its definition, graded involution 
satisfies 
\begin{equation}
\widehat{AB} = \hat{A}\hat{B} . 
\end{equation}
Consequently, conjugation satisfies 
\begin{equation}
\overline{AB} = \bar{B}\bar{A} .
\end{equation}
Table~\ref{tab:signchanges} summarizes the sign changes resulting from these 
operations on $k$-vectors: 
\begin{table}
\begin{center}
\begin{tabular}{|c|| c|c|c||}
\hline 
$k$ & $\hat{\;}$ & $\tilde{\;}$ & $\bar{\;\;}$  \\ \hline \hline
$0$ & $+$ & $+$ & $+$ \\
$1$ & $-$ & $+$ & $-$ \\
$2$ & $+$ & $-$ & $-$ \\
$3$ & $-$ & $-$ & $+$ \\ \hline 
\end{tabular}
\end{center}
\caption{Sign changes of operators.}
\label{tab:signchanges}
\end{table}
Note that paravectors can be characterized as the elements of 
$\mathcal{C}\ell_{3,0}$ or $\mathcal{C}\ell_{0,3}$ such that
\begin{equation}
\bfsf{P} = \tilde{\bfsf{P}} . 
\end{equation}

Elements of $\mathcal{C}\ell_{3,0}$ or $\mathcal{C}\ell_{0,3}$ are called
even if $A = \hat{A}$ and odd if $A = -\hat{A}$.
The set of even elements are closed under the
Clifford product, and it is called the even
subalgebra, denoted by $\mathcal{C}\ell_{3,0}^+$ or $\mathcal{C}\ell_{0,3}^+$. 
An important
characteristic of Clifford algebras is that the even subalgebras
are also Clifford algebras, and in the present case we have, 
\begin{equation}
\mathcal{C}\ell_{3,0}^+ \simeq \mathcal{C}\ell_{0,3}^+ \simeq \mathbb{H} . 
\end{equation}

The Clifford algebras $\mathcal{C}\ell_{3,0}$ and 
$\mathcal{C}\ell_{0,3}$ are two different models 
for the description of some geometric properties
of the 3D Euclidean metric space. The operations of
reflection and rotation have a simple description 
in terms of Clifford algebras. The reflection 
$\mathcal{S}_\mathbf{u}(\mathbf{v})$ 
of a vector $\mathbf{v}$ on a plane with normal vector 
$\mathbf{u}$ is given by 
\begin{equation}
\label{reflection.trans}
\mathcal{S}_\mathbf{u}(\mathbf{v}) = - \mathbf{u}\mathbf{v}\mathbf{u}^{-1} . 
\end{equation}
This expression holds for both $\mathcal{C}\ell_{3,0}$ and
$\mathcal{C}\ell_{0,3}$, the difference between these
algebras being in the sign of $\mathbf{u}^2$ in 
$\mathbf{u}^{-1} = \mathbf{u} \mathbf{u}^{-2}$. 
The expression for a rotation follows from the Cartan-Dioudonn\'e 
theorem~\cite{VazRocha}, which gives 
\begin{equation}
\label{rotation.trans}
\mathcal{R}(\mathbf{v}) = R \mathbf{u} \tilde{R} , 
\end{equation}
where $R$ is an even element such that $R \tilde{R} = \tilde{R} R = 1$. 
If $\mathbf{u}_1$ and $\mathbf{u}_2$ are unitary and orthogonal vectors 
that span a plane $\Pi$, then the above expression with 
\begin{equation}
\label{gen.rotor}
R = {\mbox e}^{\theta \mathbf{u}_1\mathbf{u}_2/2} = 
\cos{\frac{\theta}{2}} + \mathbf{u}_1\mathbf{u}_2 \sin{\frac{\theta}{2}} 
\end{equation}
is a rotation of $\mathbf{v}$ in the plane $\Pi$ by an 
angle $\theta$. The set of even elements such that $R\tilde{R} = 
\tilde{R} R = 1$ has the properties of a group. This group is
the Spin group $\operatorname{Spin}(3,0)$, or $\operatorname{Spin}(0,3)$,  
and because $\mathcal{C}\ell_{3,0}^+ \simeq \mathcal{C}\ell_{0,3}^+ 
\simeq \mathbb{H}$ we have 
\begin{equation}
\operatorname{Spin}(3,0) \simeq \operatorname{Spin}(0,3) 
\simeq S^3 , 
\end{equation}
where $S^3$ denotes the group of unit quaternions, and 
$S^3 \simeq \operatorname{SU}(2)$. Therefore, rotations
have an equivalent description in terms of the 
algebras $\mathcal{C}\ell_{3,0}$ and $\mathcal{C}\ell_{0,3}$. 
The only difference between the two cases is that, 
while in $\mathcal{C}\ell_{3,0}$ we have $\mathbf{u}\cdot 
\mathbf{v} = g(\vec{u},\vec{v})$, in $\mathcal{C}\ell_{0,3}$ 
we have $\mathbf{u}\cdot \mathbf{v} = -g(\vec{u},\vec{v})$, 
so that the rotation described by $R$ in eq.\eqref{gen.rotor} has
opposite directions in $\mathcal{C}\ell_{3,0}$ and $\mathcal{C}\ell_{0,3}$; 
in other words,  
to describe the same rotation, if we use
$\theta$ in the expression in eq.\eqref{gen.rotor} 
for $\mathcal{C}\ell_{3,0}$, we have to use 
$-\theta$  in eq.\eqref{gen.rotor} for $\mathcal{C}\ell_{0,3}$.

If rotations (and reflections) are the only geometric 
transformation we are interested, then there is absolutely 
no reason to choose $\mathcal{C}\ell_{3,0}$ over $\mathcal{C}\ell_{0,3}$ 
or vice-versa as our algebraic model for the multiplication of vectors. 
But these are not the only geometric transformations used in
computer graphics.  We will see in section~\ref{Trans.points} 
that  $\mathcal{C}\ell_{3,0}$ or $\mathcal{C}\ell_{0,3}$
are not capable of handling these other transformations.

\section{3D Euclidean space and the Clifford algebra $\mathcal{C}\ell_{3,3}$}
\label{sec.3}

Considering that the dual vector space is indissoluble from a vector space, 
let us consider it from the beginning in our algebraic model, 
and see if we can study  geometric transformations other 
than reflection and rotation within this model. 

We will take as our algebraic model the Clifford algebra $\mathcal{C}\ell_{3,3}$. 
Let us denote the generators of $\mathcal{C}\ell_{3,3}$ by $\mathbf{e}_i^+$ and $\mathbf{e}_i^-$ ($i=1,2,3$) 
subject to the following relations: 
\begin{alignat}{1}
\label{clif.33.1}
& \mathbf{e}_i^+ \mathbf{e}_j^+ + \mathbf{e}_j^+ \mathbf{e}_i^+ = 2 g_{ij} , \\
\label{clif.33.2}
& \mathbf{e}_i^- \mathbf{e}_j^- + \mathbf{e}_j^- \mathbf{e}_i^- = -2 g_{ij} , \\
\label{clif.33.3}
& \mathbf{e}_i^+ \mathbf{e}_j^- + \mathbf{e}_j^- \mathbf{e}_i^+ = 0 ,
\end{alignat}
where $g_{ij} = g(\vec{e}_i,\vec{e}_j)$. {Let us suppose that 
$g_{ij} = \delta_{ij}$.} 

The interpretation of $\mathbf{e}_i^+$ and $\mathbf{e}_i^-$ comes
with the following definition. Given vectors $\vec{e}_i$ and
covectors $\vec{e}_i{}^\ast$, we define $\gamma: V_3 \rightarrow 
\mathcal{C}\ell_{3,3}$ and $\gamma: V_3^\ast \rightarrow 
\mathcal{C}\ell_{3,3}$ as 
\begin{alignat}{1}
\label{gamma.vec}
& \gamma(\vec{e}_i) = \mathbf{e}_i = \frac{1}{2}(\mathbf{e}_i^+ + \mathbf{e}_i^-) , \\
\label{gamma.covec}
& \gamma(\vec{e}_i{}^\ast) = \mathbf{e}_i^\ast = \frac{1}{2}(\mathbf{e}_i^+ - 
\mathbf{e}_i^-) . 
\end{alignat}
Therefore we have 
\begin{equation}
\mathbf{e}_i^+ = \frac{1}{2}(\mathbf{e}_i + \mathbf{e}_i^\ast) , \qquad 
\mathbf{e}_i^- = \frac{1}{2}(\mathbf{e}_i - \mathbf{e}_i^\ast) . 
\end{equation}
The commutation relation of $\mathbf{e}_i$ and $\mathbf{e}_i^\ast$
follows from the commutation relation of
$\mathbf{e}_i^+$ and $\mathbf{e}_i^-$:
\begin{alignat}{1}
& \mathbf{e}_i \mathbf{e}_j + \mathbf{e}_j \mathbf{e}_i = 0 , \\
& \mathbf{e}_i^\ast \mathbf{e}_j^\ast + \mathbf{e}_j^\ast \mathbf{e}_i^\ast = 0 , \\
& \mathbf{e}_i \mathbf{e}_j^\ast + \mathbf{e}_j^\ast \mathbf{e}_i = g_{ij} . 
\end{alignat} 

A vector $\vec{v} = v^i \vec{e}_i$ is represented in the Clifford algebra
$\mathcal{C}\ell_{3,3}$ as 
\begin{equation}
\mathbf{v} = v^i \mathbf{e}_i . 
\end{equation}
In an obvious notation, we can also write it as 
\begin{equation}
\mathbf{v} = \frac{1}{2}(\mathbf{v}^+ + \mathbf{v}^-) = \frac{1}{2}\sum_\sigma \mathbf{v}^{\sigma_{\mathbf{v}}} ,  
\end{equation}
where $\sigma_{(\cdot)} = \{+,-\}$. 
Note that, while $(\mathbf{v}^+)^2 = |\vec{v}|^2$ and 
$(\mathbf{v}^-)^2 = -|\vec{v}|^2$, we have
\begin{equation}
(\mathbf{v})^2 = 0 . 
\end{equation}
However, we have
\begin{equation}
\mathbf{v}\mathbf{v}^\ast + \mathbf{v}^\ast \mathbf{v} = |\vec{v}|^2 , 
\end{equation}
where $\mathbf{v}^\ast$ is the representation in $\mathcal{C}\ell_{3,3}$ 
of the covector $\vec{v}^\ast$,
\begin{equation}
\mathbf{v}^\ast = v^i \mathbf{e}_i^\ast = v_i \boldsymbol{\theta}^i , 
\end{equation}
where we defined 
\begin{equation}
\boldsymbol{\theta}^i = g^{ij}\mathbf{e}_j^\ast , \qquad 
v_i = g_{ij} v^j , 
\end{equation}
where $g^{ij}$ is such that $g^{ij}g_{jk} = \delta^i_k$. 
We also have 
\begin{equation}
\mathbf{v}^\ast = \frac{1}{2}(\mathbf{v}^+ -\mathbf{v}^-) = 
\frac{1}{2}\sum_\sigma \sigma_{\mathbf{v}} 
\mathbf{v}^{\sigma_{\mathbf{v}}} . 
\end{equation}

The representation of bivectors and trivectors in this model has
to be constructed from the products of vectors $\mathbf{v}$. Then, 
unlike the case of vectors, where we have positive ($\mathbf{v}^+$) 
and negative ($\mathbf{v}^-$) parts, for a bivector we have 
positive, negative and mixed parts, that is, 
\begin{equation}
\mathbf{u}\wedge\mathbf{v} = \frac{1}{4}(\mathbf{u}^+\wedge 
\mathbf{v}^+ + \mathbf{u}^-\wedge\mathbf{v}^- + 
\mathbf{u}^+\wedge\mathbf{v}^- + \mathbf{u}^-\wedge\mathbf{v}^+) = 
\frac{1}{2^2}\sum_\sigma \mathbf{u}^{\sigma_{\mathbf{u}}}\wedge
\mathbf{v}^{\sigma_{\mathbf{v}}} . 
\end{equation}
The same happens with trivectors, that is, 
\begin{equation}
\mathbf{u}\wedge\mathbf{v}\wedge\mathbf{w} = 
\frac{1}{2^3}\sum_{\sigma} \mathbf{u}^{\sigma_{\mathbf{u}}}\wedge
\mathbf{v}^{\sigma_{\mathbf{v}}}\wedge \mathbf{w}^{\sigma_{\mathbf{w}}} .
\end{equation} 

We also observe that there is an easy algebraic way to implement the transformation 
$\mathbf{v}\mapsto \mathbf{v}^\ast$. Let $\mathbf{I}^\pm$ be 
defined as 
\begin{equation}
\mathbf{I}^+ = \mathbf{e}_1^+ \mathbf{e}_2^+ \mathbf{e}_3^+ , \qquad 
\mathbf{I}^- = \mathbf{e}_1^- \mathbf{e}_2^- \mathbf{e}_3^- .
\end{equation}
Note that 
\begin{equation}
(\mathbf{I}^+)^2 = -1 , \qquad (\mathbf{I}^-)^2 = 1 . 
\end{equation}
Then we have 
\begin{equation}
\label{dual.alg.1}
\mathbf{I}^+ \mathbf{v}^\pm = \pm \mathbf{v}^\pm \mathbf{I}^+ , \qquad 
\mathbf{I}^- \mathbf{v}^\pm = \mp \mathbf{v}^\pm \mathbf{I}^- . 
\end{equation}
Therefore, we have 
\begin{gather}
\mathbf{I}^+ \frac{1}{2}(\mathbf{v}^+ + \mathbf{v}^-) = 
\frac{1}{2}(\mathbf{v}^+ - \mathbf{v}^-)\mathbf{I}^+ , \\
\mathbf{I}^- \frac{1}{2}(\mathbf{v}^+ + \mathbf{v}^-) = 
\frac{1}{2}(-\mathbf{v}^+ - \mathbf{v}^-)\mathbf{I}^- ,
\end{gather}
from which we conclude that 
\begin{equation}
\label{dual.alg.2}
\mathbf{I}^+ \mathbf{v}(\mathbf{I}^+)^{-1} = \mathbf{v}^\ast , \qquad 
\mathbf{I}^- \mathbf{v}(\mathbf{I}^-)^{-1} = -\mathbf{v}^\ast . 
\end{equation}
Note also that, for $\mathbf{I} = \mathbf{I}^+ \mathbf{I}^-$, we have
\begin{equation}
\mathbf{I}\mathbf{v} = - \mathbf{v}\mathbf{I} . 
\end{equation}

\subsection{The Hodge Star Duality}
\label{hodge.subsec}

In the models with the Clifford algebras $\mathcal{C}\ell_{3,0}$ or 
$\mathcal{C}\ell_{0,3}$, we use the volume elements of the underlying
vector spaces to construct the Hodge duality operator. 
However, in the model with $\mathcal{C}\ell_{3,3}$, 
the object $\mathbf{I} = \mathbf{e}_1^+ \mathbf{e}_2^+ \mathbf{e}_3^+ \mathbf{e}_1^- 
\mathbf{e}_2^- \mathbf{e}_3^-$ is the volume element of the vector 
space $V_3 \oplus V_3^\ast$, and there is 
no obvious or natural interpretation for a duality operation constructed using
the volume element of the six-dimensional vector space $V_3 \oplus V_3^\ast$. 
For this reason, we use the trivector $\Omega_V$ associated with 
the three dimensional vector space $V_3$, 
given by 
\begin{equation}
\label{omega.V}
\Omega_V = \mathbf{e}_1\mathbf{e}_2\mathbf{e}_3 = 
\frac{1}{2^3 }\sum_{\sigma} \mathbf{e}_1^{\sigma_1} 
\mathbf{e}_2^{\sigma_2}\mathbf{e}_3^{\sigma_3} . 
\end{equation}
We define the Hodge star 
operator as 
\begin{equation}
\begin{split}
& \star 1 = \Omega_V , \\
& \star \mathbf{v}^{\sigma_{\mathbf{v}}} =  2 \sigma_{\mathbf{v}} \, 
\mathbf{v}^{\sigma_{\mathbf{v}}}  \cdot \Omega_V , \\
& \star (\mathbf{u}^{\sigma_{\mathbf{u}}}  \wedge \mathbf{v}^{\sigma_{\mathbf{v}}} ) = 2^2 \sigma_{\mathbf{v}} \sigma_{\mathbf{u}} \, 
\mathbf{v}^{\sigma_{\mathbf{v}}} \cdot (\mathbf{u}^{\sigma_{\mathbf{u}}} \cdot \Omega_V) , \\
& \star (\mathbf{u}^{\sigma_{\mathbf{u}}}  \wedge \mathbf{v}^{\sigma_{\mathbf{v}}}  \wedge \mathbf{w}^{\sigma_{\mathbf{w}}}) = 2^3 
\sigma_{\mathbf{u}} \sigma_{\mathbf{v}} \sigma_{\mathbf{w}} \, 
\mathbf{w}^{\sigma_{\mathbf{w}}} \cdot (\mathbf{v}^{\sigma_{\mathbf{v}}} \cdot 
(\mathbf{u}^{\sigma_{\mathbf{u}}} \cdot \Omega_V))) . 
\end{split}
\end{equation}
The reason for the factor of $2$ is that the definition of $\cdot$ uses 
a factor of $2$ in the denominator and then\footnote{Alternatively 
we could have defined $\mathbf{e}_i = \frac{1}{\sqrt{2}}(\mathbf{e}_i^+ + \mathbf{e}_i^-)$ and $\mathbf{e}_i^\ast = \frac{1}{\sqrt{2}}(\mathbf{e}_i^+ - 
\mathbf{e}_i^-)$ and then we would have obtained $\mathbf{v}^\ast \cdot 
\mathbf{v} = |\vec{v}|^2$}  $\mathbf{v}^\ast \cdot \mathbf{v} = 
\frac{1}{2}(\mathbf{v} \mathbf{v}^\ast + \mathbf{v}^\ast \mathbf{v}) = \frac{1}{2}|\vec{v}|^2$. We can collect the above expressions in a single expression as 
\begin{equation}
\label{hodge.def}
\star A_k = \langle \widetilde{A}^\ast_k \Omega_V \rangle_{3-k} . 
\end{equation}

In Section~\ref{sec.cotranslation} we will also need the inverse of the Hodge star operator,
$\star^{-1}$, which in general can be shown to be
$\star^{-1}=(-1)^{k(n-k)}\star$ for grade $k$-objects; for
$n=3$, we have for $k=0,1,2,3$, $(-1)^{0(3-0)} = 1$, $(-1)^{1(3-1)} = 1$, $(-1)^{2(3-2)} = 1$, $(-1)^{3(3-3)} = 1$, and so $ \star^{-1} = \star$ in our setting.

\begin{rem}
It is interesting to note that 
\begin{equation}
\mathbf{I} = \mathbf{I}^+ \mathbf{I}^- = 2^3 \Omega_{V^\ast} \wedge\Omega_{V} , 
\end{equation}
where $\Omega_{V^\ast} = \mathbf{e}_1^\ast \mathbf{e}_2^\ast \mathbf{e}_3^\ast$. 
Therefore, in analogy with the $\mathcal{C}\ell_{3,0}$ and $\mathcal{C}\ell_{0,3}$ 
cases, we can define a duality operation (denoted by ${\scriptstyle \maltese}$) as 
\begin{equation}
{\scriptstyle \maltese} A_k = \langle \widetilde{A}_k \mathbf{I}\rangle_{6-k} .
\end{equation}
The relation between ${\scriptstyle \maltese} A_k$ and $\star A_k$ can be 
seen as follows. We have 
\begin{equation}
{\scriptstyle \maltese} A_k = 2^3 \langle \widetilde{A}_k \Omega_{V^\ast} \wedge\Omega_{V}\rangle_{6-k} , 
\end{equation}
and if $A_k$ does not contain elements from $V^\ast$, this expression reduces to 
\begin{equation}
{\scriptstyle \maltese} A_k = 2^3 \langle \widetilde{A}_k \Omega_{V^\ast}\rangle_{3-k}\wedge \Omega_V = 
2^3 (\star A_k)^\ast \wedge \Omega_V. 
\end{equation}
Although ${\scriptstyle \maltese} A_k$ may sometimes 
be useful in some calculations, 
we do not see any natural interpretation for it. 
\end{rem}

At first glance, dealing with vectors in $\mathcal{C}\ell_{3,3}$  
seems to be less simple than with $\mathcal{C}\ell_{3,0}$ or 
$\mathcal{C}\ell_{0,3}$. For those who are used to working with 
$\mathcal{C}\ell_{3,0}$ or $\mathcal{C}\ell_{0,3}$, 
this may sound true; however, the gains from using 
$\mathcal{C}\ell_{3,3}$ are many, as we will see below. 

\section{Representation of points}
\label{vect_point}

To represent points, we
define a special point $O$, the origin. Then a 
point $P$ has coordinates $(p^1,p^2,p^3)$ in relation 
to $O$, where $p^i$ are the coordinates of the vector 
$\vec{p} = p^i \vec{e}_i$. However, the structure of 
the set of points is not the same as the structure of the
set of vectors. Although we can subtract points as 
$P - O = \overrightarrow{OP} = \vec{p}$, their sum is not a 
trivial issue. We can sum points, for example, in the context of the
barycentric calculus of M\"obius~\cite{Ungar}. Given points 
$P_0,P_1,\ldots, P_n$, the idea is to associate masses 
$m_0,m_1,\ldots, m_n$ to these points, so that a point $P$
in the $n$-simplex having these points as vertices 
have coordinates $\mu_0,\mu_1,\ldots,\mu_n$ such 
that $P = \mu_0 P_0 + \mu_1 P_1 + \cdots + \mu_n P_n$, 
with $\mu_i = m_i / M$ where $M = m_0 + m_1 + \cdots + m_n$. 
Then $\mu_0 + \mu_1 + \cdots + \mu_n = 1$. 

The inconvenient aspect of the M\"obius idea is that 
we can only add point $\mu_0 P_0 + \mu_1 P_1 + \cdots + \mu_n P_n$ 
subject to the constraint $\mu_0 + \mu_1 + \cdots + \mu_n = 1$. 
We can circumvent this with the system of mass-points described, 
for example, by Goldman~\cite{Goldman}. The idea behind mass-points 
is to introduce a mass dimension to a point $P$, representing it 
as $(m,mP)$, where $m$ is its mass. Then $(m,mP)$ is a point 
in a Grassmann space~\cite{Goldman2}. The addition 
of points in Grassmann space is $(m_0,m_0 P_0) + (m_1,m_1 P_1) + 
\cdots + (m_n,m_n P_n) = (m_0 + m_1 + \cdots + m_n, m_0P_0 + m_1 P_1 + 
\cdots + m_n P_n)$. Points in affine space are defined as points
with unit mass $(1,P)$. In this system of mass-points, mass can 
be positive or negative. 

The idea of mass-points is rich, 
and we will explore a small variation of it. 
To depart from this scenario, 
we will use a different name, introduce 
the idea of weighted points, and
resort to the concept of a paravector. 
A {\em paravector} is the sum of a scalar and a vector 
\cite{VM2018,Porteous,AACA_paravectors}. Although 
this is equivalent to introducing a new 
dimension, it is in fact more than this; 
we will identify the paravectors 
within the multivector structure of a Clifford algebra, 
in such a way that we are able to identify their subspace 
using the automorphisms of the Clifford algebra. 
We will use a paravector $\mathsf{P} = w + \vec{p}$ 
to represent a point $P$. The vector $\vec{p}$ 
will be identified with $\vec{p} = \overrightarrow{OP} = P - O$, 
where $O$ is the origin. The scalar quantity $w$ 
is the weight. The weight can be positive or negative, 
and for paravectors the sign of the weight is typically interpreted as 
the orientation of the point, although we do not make particular use of this
orientation in this paper. We will say that points with $w > 0$ have 
positive orientation, and points with $w < 0$ have 
negative orientation. Points with unitary weight are
interpreted as points in an affine space, so that  
points described by the paravectors $1 + \vec{p}$ 
and $-1 + \vec{p}$ describe points with
the same location but opposite orientations. 
We call points with unitary weight {\em affine points}. 
Given a weighted point represented by $\mathsf{P}^\prime = 
w + \vec{p}{\,}^\prime$, its location is the same as the 
affine point $\mathsf{P} = 1 + \vec{p}$ when 
$\mathsf{P}^\prime = \lambda \mathsf{P}$, with $\lambda > 0$ 
to preserve orientation. Then $\lambda = w$, 
$\vec{p}{\,}^\prime = \lambda \vec{p}$, so the location 
of the weighted point $\mathsf{P}^\prime$ in relation 
to $O$ is given by the vector $w^{-1}\vec{p}{\,}^\prime$. 
If $\mathsf{P} = -1 + \vec{p}$, then $\lambda = -w$, 
so that $w < 0$, and the location of the weighted point 
is $-w^{-1}\vec{p}{\,}^\prime = |w|^{-1}\vec{p}{\,}^\prime$. 

When we consider the Clifford algebra $\mathcal{C}\ell_{3,3}$, 
the paravector $\mathsf{P} = w + \vec{p}$ is represented 
by $\bfsf{P} = 1 + \mathbf{p}$, 
which we can also write as 
\begin{equation}
\bfsf{P} = \frac{1}{2}(\bfsf{P}^+ + \bfsf{P}^-) , 
\end{equation}
where 
\begin{equation}
\bfsf{P}^+ = 1 + \mathbf{p}^+ , \qquad 
\bfsf{P}^- = 1 + \mathbf{p}^- . 
\end{equation}

\section{Transformations of points}
\label{Trans.points}
 
We are interested in studying invertible transformations of points. Since points
are described in this formalism by means of paravectors 
$\bfsf{P} = 1 + \mathbf{p}$, we will consider algebraic 
operations of the form $U \bfsf{P} V$ whose result is 
another paravector. From the fact that paravectors $\bfsf{P}$ 
are elements such that $\bfsf{P} = \widetilde{\bfsf{P}}$, 
we have that 
\begin{equation}
U\bfsf{P} V = \tilde{V}\bfsf{P} \tilde{U} . 
\end{equation}
Since the transformation $U \bfsf{P} V$ should
be invertible, we have 
\begin{equation}
U^{-1}\tilde{V}\bfsf{P} \tilde{U}V^{-1} = \bfsf{P} , 
\end{equation}
from which we conclude, from the arbitrariness of $\bfsf{P}$, that 
\begin{equation}
V = \epsilon \tilde{U} , \qquad \epsilon = \pm 1. 
\end{equation}
Therefore, the transformation of points has to be of the form 
\begin{equation}
\label{trans.eq}
\bfsf{P}^\prime = \epsilon U \bfsf{P}\tilde{U} .
\end{equation} 
However, in $\mathcal{C}\ell_{3,3}$ 
paravectors are not the only elements that satisfy $\phi = \tilde{\phi}$. 
Elements of $\bigwedge_4\oplus \bigwedge_5$ also satisfy this 
property. So, we have to add the condition 
\begin{equation}
\label{eq.gar}
\langle U\bfsf{P}\tilde{U} \rangle_{4\oplus 5} = 
\langle U \bfsf{P}\tilde{U}\rangle_4 + 
\langle U \bfsf{P} \tilde{U}\rangle_5 = 0 
\end{equation}
to assure that the result of the transformation $\epsilon U\bfsf{P}\tilde{U}$ 
is a paravector. 

From eq.\eqref{trans.eq} we see that the weight of  $P^\prime$ is 
\begin{equation}
\langle \bfsf{P}^\prime \rangle_0 = \epsilon \langle U \tilde{U}\rangle_0 + 
\langle U \mathbf{p}\tilde{U}\rangle_0 , 
\end{equation} 
while the vector part of $P^\prime$ is
\begin{equation}
\langle \bfsf{P}^\prime \rangle_1 = \epsilon \langle U \tilde{U}\rangle_1 + 
\langle U \mathbf{p}\tilde{U}\rangle_1 , 
\end{equation}
that is, an arbitrary general transformation that mixes the 
scalar and vector parts. However, if $U$ is an even element, there 
is no such mixing. In fact, if $U$ is even, then (i) $U \mathbf{p}\tilde{U}$ 
is an odd element, so its scalar part vanishes, and (ii) $U \tilde{U}$ 
is an even element, so its vector part vanishes.

\subsection{Reflection}

\begin{thm}
\label{thm.reflection}
Let $P$ be a three dimensional point described by the
paravector $\bfsf{P}=1+\mathbf{p}$. 
The point ${P}_{\scriptscriptstyle R}$ obtained by the reflection of
${P}$  
on a plane with a unitary vector $\vec{n} = n^i \vec{e}_i$ is
given by 
\begin{equation}
\bfsf{P}_{\scriptscriptstyle R} = - N \bfsf{P}\tilde{N} , 
\end{equation}
where 
\begin{equation}
\label{refl.operator}
N = N(\vec{n}) = \mathbf{n}^+ \mathbf{n}^- ,  
\end{equation}
with 
\begin{equation}
\mathbf{n}^+ = n^i \mathbf{e}_i^+ , \quad \mathbf{n}^- = n^i \mathbf{e}_i^- . 
\end{equation}
\end{thm}
 
\medskip

\begin{proof}
 We have that 
\begin{equation}
\begin{split}
-N \bfsf{P}^+ \tilde{N} & = - \mathbf{n}^+ \mathbf{n}^- (1+ \mathbf{p}^+) 
\mathbf{n}^- \mathbf{n}^+ 
 = -\mathbf{n}^+ \mathbf{n}^- \mathbf{n}^- \mathbf{n}^+ - 
\mathbf{n}^+ \mathbf{n}^-  \mathbf{p}^+ \mathbf{n}^- \mathbf{n}^+ \\
& = \mathbf{n}^+ \mathbf{n}^+ - \mathbf{n}^+ \mathbf{p}^+ \mathbf{n}^+ = 
1 + \mathcal{S}_{\mathbf{n}^+}(\mathbf{p}^+) = 1 + \mathbf{p}^+_{\scriptscriptstyle R},
\end{split} 
\end{equation}
where $\mathbf{p}^+_{\scriptscriptstyle R} = 
\mathcal{S}_{\mathbf{n}^+}(\mathbf{p}^+)$ is the reflection 
transformation as described in eq.\eqref{reflection.trans} 
in terms of $\mathcal{C}\ell_{3,0}$. 
Analogously, 
\begin{equation}
\begin{split}
-N \bfsf{P}^- \tilde{N} & = - \mathbf{n}^+ \mathbf{n}^- (1+ \mathbf{p}^-) 
\mathbf{n}^- \mathbf{n}^+ = - \mathbf{n}^- \mathbf{n}^+ (1+ \mathbf{p}^-) 
\mathbf{n}^+ \mathbf{n}^- \\
& = -\mathbf{n}^- \mathbf{n}^- + \mathbf{n}^- \mathbf{p}^- \mathbf{n}^- = 
1 + \mathcal{S}_{\mathbf{n}^-}(\mathbf{p}^-) = 1 + 
\mathbf{p}^-_{\scriptscriptstyle R}, 
\end{split}
\end{equation}
where $\mathbf{p}^-_{\scriptscriptstyle R} = 
\mathcal{S}_{\mathbf{n}^-}(\mathbf{p}^-)$ is the reflection 
transformation as described in eq.\eqref{reflection.trans} 
in terms of $\mathcal{C}\ell_{0,3}$. When we add 
$-N \bfsf{P}^+ \tilde{N}$ and $-N \bfsf{P}^- \tilde{N}$
we obtain 
\begin{equation}
-N \bfsf{P} \tilde{N} = \frac{1}{2}(1 + \mathbf{p}^+_{\scriptscriptstyle R} + 
1 + \mathbf{p}^-_{\scriptscriptstyle R}) = 1 + \mathbf{p}_{\scriptscriptstyle R} , 
\end{equation}
which gives the reflected point $P_{\scriptscriptstyle R}$. 
\end{proof}

\medskip
\begin{rem}
In the above proof, 
the reflection operation is performed \textit{independently} in 
the positive ($\bfsf{P}^+$) and negative ($\bfsf{P}^-$) parts
of $\bfsf{P}$. This is why reflection can be equally described in 
terms of $\mathcal{C}\ell_{3,0}$ or $\mathcal{C}\ell_{0,3}$. 
Moreover, since the
composition of two reflections is a rotation (Cartan-Dieudonn\'e theorem),
this independence on the positive and negative parts of $\bfsf{P}$
is expected to happen with the rotation operation, as we
will see in section~\ref{sec:rotation}.
\end{rem}

\begin{rem}
In what follows we will consider four different transformations of points 
(theorems~\ref{thm.4},~\ref{thm.5},~\ref{thm:shear}~and~\ref{thm:scale}) 
for which we have $U \tilde{U} = 1$. Therefore
the weight of the point $\bfsf{P}$ is not changed in these
transformations, so we need
to consider only the action of the transformation of its vector part 
to study the action of the transformation on a point. 
\end{rem}

\subsection{Rotation}
\label{sec:rotation}
\begin{thm}
\label{thm.4}
Let $\mathbf{p}$, $\mathbf{u}$ and $\mathbf{v}$ represent three dimensional vectors 
with $g(\vec{u},\vec{v}) = 0$ and 
$|\vec{u}| = |\vec{v}| = 1$. Then the 
transformation $R \mathbf{p}\tilde{R}$ with 
\begin{equation}
\label{rot.operator}
R = R(\vec{u},\vec{v};\theta) = 
{\mbox e}^{\theta (\mathbf{u}^+ \mathbf{v}^+ - \mathbf{u}^- \mathbf{v}^-)/2} 
\end{equation}
is a rotation of the vector $\mathbf{p}$ by an angle $\theta$ in the
plane of $\vec{u}$ and $\vec{v}$.  
\end{thm}

\medskip

\begin{proof}
We note that 
\begin{equation}
R = {\mbox e}^{\theta (\mathbf{u}^+ \mathbf{v}^+ - \mathbf{u}^- \mathbf{v}^-)/2}
= {\mbox e}^{\theta \mathbf{u}^+ \mathbf{v}^+/2} 
{\mbox e}^{-\theta \mathbf{u}^- \mathbf{v}^-/2} =  
{\mbox e}^{-\theta \mathbf{u}^- \mathbf{v}^-/2} 
{\mbox e}^{\theta \mathbf{u}^+ \mathbf{v}^+/2} . 
\end{equation}
Then we have 
\begin{equation}
\begin{split}
R(\mathbf{p}^+ + \mathbf{p}^-)\tilde{R} & = 
{\mbox e}^{-\theta \mathbf{u}^- \mathbf{v}^-/2} 
{\mbox e}^{\theta \mathbf{u}^+ \mathbf{v}^+/2} \mathbf{p}^+ 
{\mbox e}^{+\theta \mathbf{u}^- \mathbf{v}^-/2} 
{\mbox e}^{-\theta \mathbf{u}^+ \mathbf{v}^+/2} \\
& \quad + 
{\mbox e}^{-\theta \mathbf{u}^- \mathbf{v}^-/2} 
{\mbox e}^{\theta \mathbf{u}^+ \mathbf{v}^+/2} \mathbf{p}^- 
{\mbox e}^{+\theta \mathbf{u}^- \mathbf{v}^-/2} 
{\mbox e}^{-\theta \mathbf{u}^+ \mathbf{v}^+/2} \\
& = {\mbox e}^{\theta \mathbf{u}^+ \mathbf{v}^+/2} \mathbf{p}^+ 
{\mbox e}^{-\theta \mathbf{u}^+ \mathbf{v}^+/2} + 
{\mbox e}^{-\theta \mathbf{u}^- \mathbf{v}^-/2}  \mathbf{p}^- 
{\mbox e}^{\theta \mathbf{u}^- \mathbf{v}^-/2} \\
& = 
\mathcal{R}(\mathbf{p}^+) + \mathcal{R}(\mathbf{p}^-) , 
\end{split}
\end{equation}
giving the rotation of $\mathbf{p}$. Note that, as expected, the
rotation operation acts independently on $\mathbf{p}^+$ and 
$\mathbf{p}^-$.
\end{proof}

\subsection{Hyperbolic rotation}

\begin{thm}
\label{thm.5}
Let $\mathbf{p}$, $\mathbf{u}$ and $\mathbf{v}$ represent three dimensional vectors 
with $g(\vec{u},\vec{v}) = 0$ and 
$|\vec{u}| = |\vec{v}| = 1$. Then the 
transformation $H \mathbf{p}\tilde{H}$ with 
\begin{equation}
\label{hyperb.operator}
H = H(\vec{u},\vec{v};\eta) = 
{\mbox e}^{\eta (\mathbf{u}^- \mathbf{v}^+ + \mathbf{v}^- \mathbf{u}^+)/2} 
\end{equation}
is a hyperbolic rotation of the vector $\mathbf{p}$ by an angle $\eta$ in the
plane of $\vec{u}$ and $\vec{v}$.  
\end{thm}

\medskip

\begin{proof}
 Note that 
\begin{equation}
{\mbox e}^{\eta (\mathbf{u}^- \mathbf{v}^+ + \mathbf{v}^- \mathbf{u}^+)/2} = 
{\mbox e}^{\eta \mathbf{u}^- \mathbf{v}^+/2} 
{\mbox e}^{\eta \mathbf{v}^- \mathbf{u}^+/2} .
\end{equation}
We denote 
\begin{equation}
\mathbf{p}^\pm = \mathbf{p}^\pm_\parallel + \mathbf{p}^\pm_\perp , 
\end{equation}
with 
\begin{equation}
\mathbf{p}^\pm_\parallel = p_u \mathbf{u}^\pm + p_v \mathbf{v}^\pm . 
\end{equation}
Then standard calculations give 
\begin{equation}
\begin{split}
& H \mathbf{p}^+ \tilde{H} = p_u(\cosh{\eta}\mathbf{u}^+ + \sinh\eta \mathbf{v}^-) 
+ p_v(\cosh{\eta}\mathbf{v}^+ + \sinh\eta \mathbf{u}^-) + \mathbf{p}^+_\perp ,\\
& H \mathbf{p}^- \tilde{H} = p_u(\cosh{\eta}\mathbf{u}^- + \sinh\eta \mathbf{v}^+) 
+ p_v(\cosh{\eta}\mathbf{v}^- + \sinh\eta \mathbf{u}^+) + \mathbf{p}^-_\perp , 
\end{split}
\end{equation}
and then 
\begin{equation}
H \mathbf{p}\tilde{H} = \mathbf{u} (p_u\cosh\eta + p_v\sinh\eta) + 
\mathbf{v}(p_v \cosh\eta + p_u \sinh\eta) + \mathbf{p}_\perp , 
\end{equation}
which is a hyperbolic rotation by an angle $\eta$. 
\end{proof} 
 
\medskip

\begin{rem}
The proof of the above theorem clearly shows that the 
transformation $H \mathbf{p}\tilde{H}$ mixes the 
positive ($\mathbf{p}^+$) and negative ($\mathbf{p}^-$) 
parts of $\mathbf{p}$. Therefore, unlike reflections and 
rotations, hyperbolic rotation cannot be implemented in 
$\mathcal{C}\ell_{3,0}$ or $\mathcal{C}\ell_{0,3}$ by 
a transformation of the form $H \mathbf{p}\tilde{H}$ with 
$H$ an even element of the algebra. The same mixing happens with 
the shear transformation, the non-uniform scale 
transformation, and translation, as we will see in the next three sections. 
\end{rem}

\subsection{Shear transformation}

\begin{thm}
  \label{thm:shear}
Let $\mathbf{p}$, $\mathbf{u}$ and $\mathbf{v}$ represent three dimensional vectors 
such that $g(\vec{u},\vec{v}) = 0$. Then the
transformation 
\begin{equation}
\mathbf{p}^\prime = S \mathbf{p}\tilde{S}
\end{equation}
with
\begin{equation}
\label{shear.operator}
S = S(\vec{u},\vec{v};t) = 
{\mbox e}^{t(\mathbf{u}^+ + \mathbf{u}^-)(\mathbf{v}^+ - \mathbf{v}^-)/4}
\end{equation}
is a shear in the plane spanned by $\vec{u}$ and $\vec{v}$, 
that is, 
\begin{equation}
\label{shear.eq}
S \mathbf{p}\tilde{S} = 
\mathbf{p} +  t p_v \mathbf{u} . 
\end{equation}
\end{thm}

\medskip

\begin{proof}
 When $g(\vec{u},\vec{v}) = 0$ 
we have 
\begin{equation}
[(\mathbf{u}^+ + \mathbf{u}^-)(\mathbf{v}^+ - \mathbf{v}^-)]^ 2 = 0 , 
\end{equation}
and then 
\begin{equation}
{\mbox e}^{t(\mathbf{u}^+ + \mathbf{u}^-)(\mathbf{v}^+ - \mathbf{v}^-)/4} = 
1 + \frac{t}{4}(\mathbf{u}^+ + \mathbf{u}^-)(\mathbf{v}^+ - \mathbf{v}^-) . 
\end{equation}
Then the usual calculations give 
\begin{equation}
\begin{split}
& S\mathbf{p}^+ \tilde{S} = \mathbf{p}^+ + 
\frac{t}{2} p_v(\mathbf{u}^+ + \mathbf{v}^+) -\frac{t}{2}p_u(\mathbf{v}^+ - 
\mathbf{v}^-) , \\
& S\mathbf{p}^- \tilde{S} = \mathbf{p}^- + 
\frac{t}{2} p_v(\mathbf{u}^+ + \mathbf{v}^+) + \frac{t}{2}p_u(\mathbf{v}^+ - 
\mathbf{v}^-) , 
\end{split}
\end{equation}
and adding the expressions we obtain eq.\eqref{shear.eq}. 
\end{proof}

\subsection{Non-uniform scale transformations}

\begin{thm}
\label{thm:scale}
Let $\mathbf{p}$ and $\mathbf{u}$ represent three dimensional vectors, 
and suppose $\vec{u}$ is unitary. Then the
transformation 
\begin{equation}
\mathbf{p}^\prime = D \mathbf{p}\tilde{D}
\end{equation}
with
\begin{equation}
\label{non-u.operator}
D = D(\vec{u};t) = {\mbox e}^{t\mathbf{u}^-\mathbf{u}^+/2}
\end{equation}
is a non-uniform scale transformation of the
component $\mathbf{p}_\parallel$ of $\mathbf{p}$ in the direction of $\mathbf{u}$, that
is, 
\begin{equation}
\label{scale.eq}
D \mathbf{p}\tilde{D} = \mathbf{p}_\perp + 
{\mbox e}^{t} \mathbf{p}_\parallel .
\end{equation}
\end{thm}

\medskip

\begin{proof}
 We have $(\mathbf{u}^-\mathbf{u}^+)^2 = 1$. 
Then, if we denote by $\mathbf{p}_\parallel^\pm$ the 
component of $\mathbf{p}^\pm$ in the direction of $\mathbf{u}^\pm$, and 
by $\mathbf{p}_\perp^\pm$ the 
component of $\mathbf{p}^\pm$ orthogonal to $\mathbf{u}^\pm$, 
we obtain, after the usual calculations, that 
\begin{equation}
\begin{split}
& D \mathbf{p}^+ \tilde{D} = \cosh{t}\mathbf{p}_\parallel^+ 
+ g(\vec{p},\vec{u}) \sinh{t} \mathbf{u}^- + \mathbf{p}^+_\perp , \\
& D \mathbf{p}^- \tilde{D} = \cosh{t}\mathbf{p}_\parallel^- 
+ g(\vec{p},\vec{u}) \sinh{t} \mathbf{u}^+ + \mathbf{p}^-_\perp .
\end{split}
\end{equation}
When we add these expression we obtain 
\begin{equation}
D\mathbf{p} \tilde{D} = (\cosh{t}+\sinh{t})\mathbf{p}_\parallel + 
\mathbf{p}_\perp , 
\end{equation}
which is eq.\eqref{scale.eq}. 
\end{proof}

\subsection{Translation}

All transformations we have studied so far involve 
products of $\mathbf{u}^+$ and $\mathbf{v}^-$ 
and have no action on the scalar part of $\bfsf{P}$. 
Let us see now another possibility. 

\begin{thm}
Let $\bfsf{P}= 1 +\mathbf{p}$ represent a three dimensional point 
and $\mathbf{v} = \frac{1}{2}(\mathbf{v}^+ + \mathbf{v}^-)$ represents a three dimensional vector. 
Then the
transformation 
\begin{equation}
\bfsf{P}^\prime = T \bfsf{P}\tilde{T}
\end{equation}
with $T$ given by 
\begin{equation}
\label{trans.translation}
T = T(\vec{v}) = {\mbox e}^{\mathbf{v}/2} 
\end{equation} 
is a translation of the point $\bfsf{P}$ by the
vector $\mathbf{v}$, 
that is, 
\begin{equation}
T \bfsf{P}\tilde{T} = 
\bfsf{P} + \mathbf{v} . 
\end{equation}
\end{thm}

\medskip

\begin{proof}
 Since $\mathbf{v}^2 = 0$ we have 
\begin{equation}
{\mbox e}^{\mathbf{v}/2} = 1 + \frac{1}{2}\mathbf{v} . 
\end{equation}
The action of $T$ on the positive and negative vector parts of $\bfsf{P}$ are 
\begin{equation}
\begin{split}
& {\mbox e}^{\mathbf{v}/2}\mathbf{p}^+ {\mbox e}^{\mathbf{v}/2} = 
\mathbf{p}^+ + g(\vec{p},\vec{v}) + \frac{1}{2}g(\vec{p},\vec{v})\mathbf{v} , \\
& {\mbox e}^{\mathbf{v}/2}\mathbf{p}^- {\mbox e}^{\mathbf{v}/2} = 
\mathbf{p}^- - g(\vec{p},\vec{v}) - \frac{1}{2}g(\vec{p},\vec{v})\mathbf{v} . 
\end{split}
\end{equation}
Then we have 
\begin{equation}
{\mbox e}^{\mathbf{v}/2}\mathbf{p}   {\mbox e}^{\mathbf{v}/2} = \mathbf{p} , 
\end{equation}
that is, this transformation $T$ leaves vectors unchanged. 
However, the action of $T$ on the scalar part of 
$\bfsf{P}$ is 
\begin{equation}
{\mbox e}^{\mathbf{v}/2} 1   {\mbox e}^{\mathbf{v}/2} = 
\left(1 + \frac{1}{2}\mathbf{v}\right) 1 \left(1 + \frac{1}{2}\mathbf{v}\right) = 
1 + \mathbf{v} . 
\end{equation}
Then we have 
\begin{equation}
{\mbox e}^{\mathbf{v}/2}\bfsf{P} \widetilde{{\mbox e}^{\mathbf{v}/2}} = 
\bfsf{P} + \mathbf{v} , 
\end{equation}
which is the translation of the point  $\bfsf{P}$ by the vector $\mathbf{v}$. 
\end{proof}

\begin{rem}
  As with our exterior algebra model, $\bfsf{P} = e^{\mathbf{p}/2} O e^{\mathbf{p}/2} = (e^{\mathbf{p}/2})^2$
  (since $O=1$), and $e^{\mathbf{p}/2}$ is a representation for the square root of a point.
\end{rem}
\subsection{Cotranslation}
\label{sec.cotranslation}

The cotranslation is an operation defined in \cite{VM2018}, 
where it was used to write
perspective projective in terms of 
the composition of the translation and the 
cotranslation transformations. An important difference from 
this transformation over previous ones is that cotranslation uses 
the Hodge star duality in its definition. 

\begin{thm}
Let $\bfsf{P}= 1 +\mathbf{p}$ represent a three dimensional point 
and $\mathbf{v} = \frac{1}{2}(\mathbf{v}^+ + \mathbf{v}^-)$ represent a three dimensional vector. 
Then the cotranslation 
transformation, 
\begin{equation}
\bfsf{P}^\prime = \star^{-1}[ T (\star \bfsf{P})\tilde{T}]
\end{equation}
with $T$ given by 
\begin{equation}
T = {\mbox e}^{\mathbf{v}/2} 
\end{equation} 
is such that 
\begin{equation}
\bfsf{P}^\prime = \bfsf{P} + g(\vec{p},\vec{v}) .
\end{equation}
\end{thm}

\medskip

\begin{proof}
 We first note that
\begin{equation}
\mathbf{p}^+ \cdot \mathbf{e}_i^{\sigma_i} = 
\begin{cases} p_i , \; & \sigma_i = + , \\
0 , \; & \sigma_i = - , \end{cases} \qquad 
\mathbf{p}^- \cdot \mathbf{e}_i^{\sigma_i} = 
\begin{cases} 0 , \; & \sigma_i = + , \\
-p_i , \; & \sigma_i = - , \end{cases}
\end{equation}
which gives
\begin{equation}
\label{esc.prod}
(\mathbf{p}^+ - \mathbf{p}^-)\cdot \mathbf{e}_i^{\sigma_i} 
 = \sum_{\sigma^\prime} \sigma^\prime_{\mathbf{p}} \mathbf{p}^{\sigma^\prime_{\mathbf{p}}} 
 \cdot \mathbf{e}_i^{\sigma_i} = 
p_i . 
\end{equation} 
Let us calculate 
\begin{equation}
\star \bfsf{P} = \star 1 + \star \mathbf{p} = \Omega_V + 2\mathbf{p}^\ast \cdot \Omega_V .
\end{equation}
We have 
\begin{equation}
\begin{split} 
\star \mathbf{p} & = 2 \left(\frac{1}{2} \sum_{\sigma^\prime} \sigma^\prime_{\mathbf{p}} \mathbf{p}^{\sigma^\prime_{\mathbf{p}}} \right) \cdot \left( \frac{1}{2^3} \sum_\sigma 
\mathbf{e}_1^{\sigma_1} \wedge \mathbf{e}_2^{\sigma_2} \wedge \mathbf{e}_3^{\sigma_3} 
\right) \\
& = p_1 \frac{1}{2^2} \sum_\sigma \mathbf{e}_2^{\sigma_2}\wedge \mathbf{e}_3^{\sigma_3} 
- p_2 \frac{1}{2^2} \sum_\sigma \mathbf{e}_1^{\sigma_1}\wedge \mathbf{e}_3^{\sigma_3} 
+ p_3 \frac{1}{2^2} \sum_\sigma \mathbf{e}_1^{\sigma_1}\wedge \mathbf{e}_2^{\sigma_2} , 
\end{split}
\end{equation}
where we have used the property given in eq.\eqref{Leibniz} and 
because eq.\eqref{esc.prod} is counted for $\sigma_i = +$ and 
$\sigma_i = -$ there is an additional factor 2 in the final expression. For 
$T = {\mbox e}^{\mathbf{v}/2} = 1 + \mathbf{v}/2$ we have 
\begin{equation}
T(\mathbf{e}_i^{\sigma_i}\wedge\mathbf{e}_j^{\sigma_j})\tilde{T} = 
\mathbf{e}_i^{\sigma_i}\wedge\mathbf{e}_j^{\sigma_j} + 
\mathbf{e}_i^{\sigma_i}\wedge\mathbf{e}_j^{\sigma_j}\wedge\mathbf{v} - 
\frac{1}{4}\sigma_j v_j \mathbf{e}_i^{\sigma_i}\wedge\mathbf{v} + 
\frac{1}{4}\sigma_i v_i \mathbf{e}_j^{\sigma_j}\wedge\mathbf{v} , 
\end{equation}
and then 
\begin{equation}
T\left(\sum_\sigma \mathbf{e}_i^{\sigma_i}\wedge\mathbf{e}_j^{\sigma_j}\right)\tilde{T} = 
\sum_\sigma\mathbf{e}_i^{\sigma_i}\wedge\mathbf{e}_j^{\sigma_j} + 
\sum_\sigma \mathbf{e}_i^{\sigma_i}\wedge\mathbf{e}_j^{\sigma_j}\wedge\mathbf{v} ,
\end{equation}
because the contribution of the last two terms are cancelled in pairs 
due to the different signs in $\sigma_i$ and $\sigma_j$. 
Moreover 
\begin{equation}
T \Omega_V \tilde{T} = \Omega_V , 
\end{equation}
and it follows that 
\begin{equation}
\begin{split}
T(\star \bfsf{P})\tilde{T} & = \Omega_V + 
p_1 \frac{1}{2^2}\sum_\sigma \mathbf{e}_2^{\sigma_2}\wedge 
\mathbf{e}_3^{\sigma_3}\wedge(1+\mathbf{v}) - 
p_2 \frac{1}{2^2}\sum_\sigma \mathbf{e}_1^{\sigma_1}\wedge 
\mathbf{e}_3^{\sigma_3}\wedge(1+\mathbf{v}) \\
& \quad + 
p_3 \frac{1}{2^2}\sum_\sigma \mathbf{e}_1^{\sigma_1}\wedge 
\mathbf{e}_2^{\sigma_2}\wedge(1+\mathbf{v}) . 
\end{split}
\end{equation}
To obtain $\star(T(\star \bfsf{P})\tilde{T})$, 
we first calculate 
\begin{equation}
\begin{split}
& \star\left( \sum_\sigma \mathbf{e}_i^{\sigma_i}\wedge \mathbf{e}_j^{\sigma_j}\right) = 
2^2 \sum_\sigma \sigma_j\mathbf{e}_j^{\sigma_j}\cdot 
\left[ \sigma_i\mathbf{e}_i^{\sigma_i} \cdot 
\left( \frac{1}{2^3} \sum_\sigma^\prime 
\mathbf{e}_1^{\sigma^\prime_{1}} \wedge \mathbf{e}_2^{\sigma^\prime_{2}} \wedge \mathbf{e}_3^{\sigma^\prime_{3}} 
\right) \right] \\
& = \frac{1}{2}\sum_{\sigma,\sigma^\prime}\bigg(
\delta^{\sigma,\sigma^\prime}_{[ij,23]} 
\mathbf{e}_1^{\sigma^\prime_{1}}  - 
\delta^{\sigma,\sigma^\prime}_{[ij,13]} 
\mathbf{e}_2^{\sigma^\prime_{2}} + 
\delta^{\sigma,\sigma^\prime}_{[ij,12]} 
\mathbf{e}_3^{\sigma^\prime_{3}} \bigg) ,  
\end{split}
\end{equation}
where we defined 
\begin{equation}
\delta^{\sigma,\sigma^\prime}_{[ij,ab]} = 
\delta^{\sigma_i\sigma^\prime_a}\delta^{\sigma_j\sigma^\prime_b}
\delta_{ia}\delta_{ib} - 
\delta^{\sigma_i\sigma^\prime_b}\delta^{\sigma_j\sigma^\prime_a}
\delta_{ib}\delta_{ia} .
\end{equation}
For example, 
\begin{equation}
 \star\left( \sum_\sigma \mathbf{e}_1^{\sigma_1}\wedge \mathbf{e}_2^{\sigma_2}\right) = 
 \frac{1}{2}\delta^{\sigma,\sigma^\prime}_{[12,12]} 
\mathbf{e}_3^{\sigma^\prime_{3}} = 2 \sum_{\sigma^\prime} \mathbf{e}_3^{\sigma^\prime_{3}} . 
\end{equation}
Moreover, we have 
\begin{equation}
\star(\mathbf{u}\wedge\mathbf{v}\wedge\mathbf{w}) = 2\mathbf{w}^\ast \cdot 
(\star(\mathbf{u}\wedge\mathbf{v})) , 
\end{equation}
and then 
\begin{equation}
\begin{split}
& \star\left( \sum_\sigma \mathbf{e}_i^{\sigma_i}\wedge \mathbf{e}_j^{\sigma_j} 
\wedge \mathbf{v} \right) = 2\mathbf{v}^\ast \cdot 
\star\left( \sum_\sigma \mathbf{e}_i^{\sigma_i}\wedge \mathbf{e}_j^{\sigma_j}\right)\\
&  = 
\sum_{\sigma^\prime}\sigma_{\mathbf{v}}\mathbf{v}^{\sigma^\prime_{\mathbf{v}}} 
\cdot \star\left( \sum_\sigma \mathbf{e}_i^{\sigma_i}\wedge \mathbf{e}_j^{\sigma_j}\right)\\ & = \sum_{\sigma,\sigma^\prime} \left( 
v_1 \delta^{\sigma,\sigma^\prime}_{[ij,23]}  - 
v_2 \delta^{\sigma,\sigma^\prime}_{[ij,13]} +
v_3  
\delta^{\sigma,\sigma^\prime}_{[ij,12]} \right) . 
\end{split}
\end{equation}
For example, 
\begin{equation}
\star\left( \sum_\sigma \mathbf{e}_1^{\sigma_1}\wedge \mathbf{e}_2^{\sigma_2} 
\wedge \mathbf{v} \right) = v_3 \sum_{\sigma,\sigma^\prime} 
\delta^{\sigma,\sigma^\prime}_{[ij,12]}  = 
4 v_3 . 
\end{equation}
Finally, 
\begin{equation}
\begin{split}
\star[T(\star\bfsf{P})\tilde{T}] & = 
\star\Omega_V + p_1\frac{1}{2^2}\left(2\sum_{\sigma^\prime}\mathbf{e}_1^{\sigma^\prime_{1}} 
+ 4 v_1\right) - 
p_2\frac{1}{2^2}\left(-2\sum_{\sigma^\prime}\mathbf{e}_2^{\sigma^\prime_{2}} 
- 4 v_2\right)\\
& \quad  + p_3\frac{1}{2^2}\left(2\sum_{\sigma^\prime}\mathbf{e}_3^{\sigma^\prime_{3}} 
+ 4 v_3\right) \\
& = 1 + \frac{1}{2}(\mathbf{p}^+ + \mathbf{p}^-) + p_1 v_1 + p_2 v_2 + p_3 v_3 \\
& = 
1 + g(\vec{p},\vec{v}) + \mathbf{p} = \bfsf{P} + g(\vec{p},\vec{v}) ,  
\end{split}
\end{equation}
which is the result. 
\end{proof}

\subsubsection{Perspective Projection Through the Composition of the Translation and the Cotranslation Transformations}

\medskip

While the previous transformations were well-known, the same does 
not happen with cotranslation. However, there is an important 
transformation that can be written using translation and 
cotranslation. It is the perspective projection, as showed 
in \cite{VM2018}. 
Consider two points $P$ and
$E$, described by the paravectors $\bfsf{P} = \mathbf{1} + \mathbf{p}$ and 
$\bfsf{E} = \mathbf{1} + \mathbf{e}$, and a plane with a normal vector $\vec{n}$ and
plane equation $\vec{x}\cdot\vec{n} = c$. 
We introduce the notation 
\begin{equation}
\mathfrak{T}_{\vec{v}}(\bfsf{P}) = {\mbox e}^{\mathbf{v}/2}
\bfsf{P}{\mbox e}^{\mathbf{v}/2} , \qquad 
\mathfrak{W}_{\vec{v}}(\bfsf{P}) = \star^{-1}[
{\mbox e}^{\mathbf{v}/2}
(\star \bfsf{P}){\mbox e}^{\mathbf{v}/2} ]. 
\end{equation}

\begin{thm}
Let $\mathbf{n}$ describe the normal to the perspective plane $\mathcal{P}$ 
with equation $\vec{x}\cdot \vec{n}  = c$, and let $\bfsf{E}$ and $\bfsf{P}$ 
be the paravectors representing the eye point $E$ and 
an arbitrary point $P$ in three dimensional space, respectively. Then 
$\bfsf{P}_0$ given by 
\begin{equation}
\bfsf{P}_0 = 
(\mathfrak{T}_{\vec{e}}\circ \mathfrak{W}_{\vec{n}/a}
\circ \mathfrak{T}_{\vec{e}}^{-1})(\bfsf{P}-\bfsf{E})
\end{equation}
where $a = c + \vec{n}\cdot\vec{e}$, 
is a weighted point in the perspective plane $\mathcal{P}$ located
at the perspective projection of $P$ from the eye point $E$ if
$P$ is located in front of $E$; if $P$ is located behind
$E$, then this same location corresponds to 
the weighted point $\bar{\bfsf{P}}_0$. 
\end{thm}

\medskip

The proof is similar to the one in \cite{VM2018}, and we omit it here.

\subsubsection{Pseudo-Perspective Projection Through Cotranslation}

A second application of cotranslation is pseudo-perspective projection.
Given an eye point $E$ looking in a direction $\vec{n}$, we wish to map the eye 
point $E$ to a point at infinity, and in particular, we want $E$ to map to 
$\pm\vec{n}$~\cite{DGM17}. In \cite{VM2018} we proved the following:

\begin{thm}
  Let $\mathbf{n}$ be a unit vector, and let $\mathbf{E}=\mathbf{1}-\mathbf{n}$.
  Then
  $\mathfrak{W}_{\vec{n}}(\mathbf{P})$ transforms the eye point $\mathbf{E}$ to the point at infinity in the direction $-\mathbf{n}$ and transforms a viewing frustum to a rectangular box.
\end{thm}

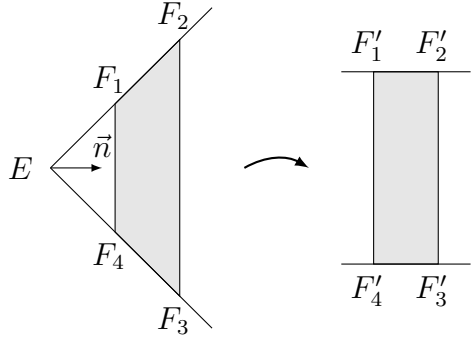
\begin{figure}
\begin{center}
\begin{tikzpicture}[>=latex,xscale = 0.85,yscale=0.85]
\fill [gray!20]
      (1, 1) -- (2,2) -- (2,-2) -- (1,-1) -- cycle;
\draw (1, 1) -- (2,2) -- (2,-2) -- (1,-1) -- cycle;
\draw (0, 0) -- (2.5,2.5);
\draw (0,0) -- (2.5,-2.5);
\node [above] at (0.9,1) {$F_1$};
\node [above] at (1.9,2) {$F_2$};
\node [below] at (1.9,-2) {$F_3$};
\node [below] at (0.9,-1) {$F_4$};
\node [left] at (-0.1,0) {$E$};
\draw [->] (0,0) -- (0.8,0);
\node [above] at (0.8,0) {$\vec{n}$} ;
\draw [thick, ->] (3,0) to [bend left] (4,0);
\fill [gray!20]
      (5, 1.5) -- (6,1.5) -- (6,-1.5) -- (5,-1.5) -- cycle;
      \draw  (5, 1.5) -- (6,1.5) -- (6,-1.5) -- (5,-1.5) -- cycle;
\draw (4.5,1.5) -- (6.5,1.5);
\draw (4.5,-1.5) -- (6.5,-1.5);
\node [above] at (4.9,1.5) {$F_1^\prime$};
\node [above] at (5.9,1.5) {$F_2^\prime$};
\node [below] at (5.9,-1.5) {$F_3^\prime$};
\node [below] at (4.9,-1.5) {$F_4^\prime$};
\end{tikzpicture}
\end{center}
  \caption{Mapping a truncated viewing pyramid to a box.  From~\cite{VM2018}.}
  \label{pseudopersp}
\end{figure}

We reproduce a figure from~\cite{VM2018} to illustrate the mapping (Figure~\ref{pseudopersp}); the proof of this result is similar to the one in~\cite{VM2018}
and the interested reader is referred to that paper for details.

\subsection{Composition of Transformations}

Rotations, hyperbolic rotations, shear transformations, 
non-uniform scale transformations and translations are 
transformations of the form 
\begin{equation}
\label{type.1}  
\bfsf{P}^\prime = U \bfsf{P} \tilde{U} , 
\end{equation}
while cotranslation is a transformation of the form 
\begin{equation}
\label{type.2}
\bfsf{P}^\prime = \star^{-1}[T (\star \bfsf{P}) \tilde{T}] .  
\end{equation}
A natural question we can ask is whether we can compose these 
transformations, that is, to write a composition of them in the 
form of a single transformation. Obviously the 
transformations of the form $\bfsf{P}^\prime = U_1 \bfsf{P} \tilde{U}_1$ 
and $\bfsf{P}^{\prime\prime} = U_2 \bfsf{P}^\prime \tilde{U}_2$ 
can be composed as $\bfsf{P}^{\prime\prime} = U_{21} \bfsf{P} \tilde{U}_{21}$ 
with $U_{21} = U_2 U_1$, as well as transformations of the 
form $\bfsf{P}^\prime = \star^{-1}[T_1 (\star \bfsf{P}) \tilde{T}_1]$ 
and $\bfsf{P}^{\prime\prime} = \star^{-1}[T_2 (\star \bfsf{P}^\prime) \tilde{T}_2]$ 
can be composed as $\bfsf{P}^{\prime\prime} = \star^{-1}[T_{21} (\star \bfsf{P}) 
\tilde{T}_{21}]$ with $T_{21} = T_2 T_1$. The question, therefore, 
is whether transformations of different types can be composed in a 
single expression. We do not expect, however, that an arbitrary 
transformation of the form of eq.\eqref{type.1} can be written in form 
of eq.\eqref{type.2}, or vice versa, because if that were possible 
one form of transformation would suffice to treat all cases. 

\begin{thm}
\label{theorem.hodge}
Let $U$ be such that:
\begin{alignat*}{2}
& \textrm{(i)} & \quad & U \tilde{U}=1 , \\
& \textrm{(ii)} & \quad & U A_k \tilde{U} = \langle U A_k 
\tilde{U}\rangle_k ,\\ 
& \textrm{(iii)} & \quad & U^\ast  A_k \tilde{U}^\ast = \langle U^\ast 
A_k \tilde{U}^\ast\rangle_k ,\\
& \textrm{(iv)} & \quad & \tilde{U}^\ast \Omega_V U^\ast = \lambda^2 \Omega_V , 
\end{alignat*}
where $A_k$ is a $k$-vector, $\Omega_V = \mathbf{e}_1\mathbf{e}_2\mathbf{e}_3$, 
$\lambda$ is a constant and $U^\ast = \mathbf{I}^+ U (\mathbf{I}^+)^{-1}$. 
Then 
\begin{equation}
\label{hodge.1}
\star (U A_k \tilde{U}) = 
U^\prime (\star A_k ) \tilde{U}^\prime , 
\end{equation} 
where $U^\prime$ is defined as 
\begin{equation}
\label{def.Uprime}
U^\prime = \lambda U^\ast. 
\end{equation}
\end{thm}

\medskip

\begin{proof}
 Using the definition in eq.\eqref{hodge.def}, 
and conditions (i)-(iv), we have 
\begin{equation}
\begin{split}
\star (U  A_k \tilde{U}) & = 
\star \langle U A_k \tilde{U}\rangle_k = 
\langle U^\ast \tilde{A}_k^\ast \tilde{U}^\ast \Omega_V \rangle_{3-k} \\
& =\langle U^\ast \tilde{A}_k^\ast \Omega_V \lambda^2 \tilde{U}^\ast\rangle_{3-k} = 
\lambda^2 U^\ast \langle \tilde{A}_k^\ast \Omega_V \rangle_{3-k} \tilde{U}^\ast , 
\end{split}
\end{equation}
from which eq.\eqref{hodge.1} follows. 
\end{proof}

\medskip 

\begin{cor}
If $U$ satisfy the conditions (i)-(iv) in Theorem~\ref{theorem.hodge}, then 
the transformation of the form $\bfsf{P}^\prime = U \bfsf{P}\tilde{U}$ 
is equivalent to the transformation of the 
form $\bfsf{P}^\prime = \star^{-1}[U^\prime (\star \bfsf{P}) \tilde{U}^\prime]$, 
with $U^\prime$ as in eq.\eqref{def.Uprime}. 
\end{cor}

\medskip

\begin{proof}
 Under the conditions for $U$, we have 
$$
\bfsf{P}^\prime = U (1 + \mathbf{p})\tilde{U} = 1 + U \mathbf{p}\tilde{U} , 
$$
where $U\mathbf{p}\tilde{U}$ is a vector. On the other hand, 
$$
U^\prime (\star \bfsf{P})\tilde{U}^\prime = U^\prime (\star 1 + 
\star \mathbf{p})\tilde{U}^\prime = 
U^\prime \Omega_V \tilde{U}^\prime + U^\prime (\star \mathbf{p})\tilde{U}^\prime . 
$$
But the result of Theorem~\ref{theorem.hodge} for $k=0$ and $k=1$ gives 
\begin{equation*}
\begin{split}
& U^\prime \Omega_V \tilde{U}^\prime = U^\prime (\star 1) \tilde{U}^\prime = 
\star (U 1 \tilde{U}) = \star 1 , \\
& U^\prime (\star p) \tilde{U}^\prime = \star (U\mathbf{p} \tilde{U}) , 
\end{split}
\end{equation*}
and then 
$$
U^\prime (\star \bfsf{P})\tilde{U}^\prime = \star(1 + U \mathbf{p}\tilde{U}) , 
$$
and comparing this with the above expression for $\bfsf{P}^\prime$, we conclude that 
\begin{equation*}
\bfsf{P}^\prime = U\bfsf{P}\tilde{U} = \star^{-1}[U^\prime (\star \bfsf{P})\tilde{U}^\prime] ,  
\end{equation*}
which is the result. 
\end{proof}

\medskip

\begin{thm}
Reflections, rotations, hyperbolic rotations, shear transformations and  
non-uniform scaling transformations can be described equivalently in terms 
of transformations of the form $U\bfsf{P}\tilde{U}$ or of the
form $\star^{-1}[U^\prime(\star \bfsf{P})\tilde{U}^\prime)]$, with 
operators $U$ and $U^\prime$ given by \smallskip 
\renewcommand{\arraystretch}{1.3}
\begin{center}
\begin{tabular}{|c|c|c||}
\hline 
 $U$ &  $U^\prime$  \\ \hline \hline
 $N(\vec{n})$ &  $-N(\vec{n})$  \\ \hline
 $R(\vec{u},\vec{v};\theta)$  & $R(\vec{u},\vec{v};\theta)$ 
 \\ \hline 
 $H(\vec{u},\vec{v};\eta)$ & $H(\vec{u},\vec{v};-\eta)$ \\ \hline 
 $S(\vec{u},\vec{v};t)$ &  $S(\vec{v},\vec{u};-t)$ \\ \hline 
 $D(\vec{u};t)$ &  ${\mbox e}^{-t/2}D(\vec{u};-t)$ \\ \hline 
\end{tabular}
\end{center}\smallskip
where $N$, $R$, $H$, $S$ and $D$ are given by eq.\eqref{refl.operator}, 
eq.\eqref{rot.operator}, eq.\eqref{hyperb.operator}, eq.\eqref{shear.operator} and 
eq.\eqref{non-u.operator}, respectively. 
\end{thm}

\medskip

\begin{proof}
 First of all, we note that 
$N$, $R$, $H$, $S$ and $D$ operators satisfy conditions (i) and (ii) 
in Theorem~\ref{theorem.hodge}. We can easily find the expression 
for $U^\ast$ using eq.\eqref{dual.alg.1}-\eqref{dual.alg.2}. It follows 
that  
\begin{gather}
[N(\vec{n})]^\ast = -N\vec{n} , \\
[R(\vec{u},\vec{v};\theta)]^\ast =  R(\vec{u},\vec{v};\theta) ,\\
[H(\vec{u},\vec{v};\eta)]^\ast = H(\vec{u},\vec{v};-\eta) ,\\  
[S(\vec{u},\vec{v};t)]^\ast = S(\vec{v},\vec{u};-t) ,\\ 
[D(\vec{u};t)]^\ast =  D(\vec{u};-t) , 
\end{gather}
and, since $U^\ast$ operators generate the same kind of transformations
as $U$ operators, condition (iii) is also satisfied. Explicit calculations 
also show that 
\begin{equation}
U^\ast \Omega_V \tilde{U}^\ast = \Omega_V ,  \;\text{for}\; 
U^\ast = \{N^\ast,R^\ast,H^\ast,S^\ast\} , 
\end{equation}
but 
\begin{equation}
D^\ast \Omega_V \tilde{D}^\ast = {\mbox e}^{-t} \Omega_V , 
\end{equation}
that is, condition (iv) is satisfied with $\lambda = 1$ for 
$U  = \{N,R,H,S\}$ and with $\lambda = {\mbox e}^{-t/2}$ for 
$U = D$. 
\end{proof} 

\medskip

\begin{rem} The translation operator does not satisfy 
the conditions of Theorem~\ref{theorem.hodge}---in fact, 
it only satisfies condition (i). However, since 
those conditions are sufficient ones, we cannot conclude
from Theorem~\ref{theorem.hodge} that there is no transformation 
of the form $\star^{-1}[U^\prime(\star \bfsf{P})\tilde{U}^\prime)]$
equivalent to a translation $T\bfsf{P}\tilde{T}$. 
However, it is not difficult to see that
there translation cannot be written in the form $\star^{-1}[U^\prime(\star \bfsf{P})\tilde{U}^\prime)]$.
For example, consider
translations of the origin along the 
direction $\mathbf{e}_1$. These points are 
$\bfsf{P} = {\mbox e}^{t\mathbf{e}_1/2} 1 {\mbox e}^{t\mathbf{e}_1/2} = 
1 + t \mathbf{e}_1$. Then $\star \bfsf{P} = \Omega_V + t\mathbf{e}_2\mathbf{e}_3$. 
If we could write translation in the form $\star^{-1}[U^\prime(\star \bfsf{P})
\tilde{U}^\prime)]$, then the operator $U$ is to be such that 
$$ 
U^\prime \star 1 \tilde{U}^\prime = U^\prime \Omega_V \tilde{U}^\prime = 
\Omega_V + t\mathbf{e}_2\mathbf{e}_3 . 
$$
To obtain the bivector $\mathbf{e}_2\mathbf{e}_3$ from 
the trivector $\Omega_V$, the operator $U^\prime$ has to contain 
a term with $\mathbf{e}_1^\ast$. Supposing $t$ sufficiently small, 
and being connect to the identity, let us write it as 
$U^\prime = 1 + t\mathbf{e}_1^\ast$. In this case 
$$
U^\prime \Omega_V \tilde{U}^\prime = \Omega_V + t \mathbf{e}_2\mathbf{e}_3 
+ t^2 \mathbf{e}_1^\ast \mathbf{e}_2\mathbf{e}_3 , 
$$
that is, we obtain an undesirable extra term $t^2 \mathbf{e}_1^\ast \mathbf{e}_2\mathbf{e}_3$. 
Then, if we try to add a term to cancel this extra term, 
this new term 
cannot contain $\mathbf{e}_1$ or $\mathbf{e}_2$ or $\mathbf{e}_3$, 
because $\mathbf{e}_i \Omega_V = \Omega_V \mathbf{e}_i = 0$ for 
$i=1,2,3$, and so it has to contain only terms involving $\mathbf{e}_i^\ast$, 
which does not cancel the extra term (but adds extra terms). Translations 
cannot, therefore, be written in the same form as cotranslations---and,
of course, cotranslations cannot be written in the same form
as translations.  

Thus, we may express reflections, rotations, shear transformations and non-uniform scaling transformations in the form $\bfsf{P}^\prime = \star^{-1}[T (\star \bfsf{P}) \tilde{T}]$ and combine them with cotranslation to obtain a single versor to apply to points.  However, to use translations with cotranslation, we are limited in which transformations we can combine.  For example, in the computer graphics setting, where we have a sequence of affine transformations followed by a single perspective transformations, we have to apply two versors to our points: a versor $A$ that is the composition of the affine transformations, and a second versor $T$ that is the cotranslation representing the perspective projection:
\begin{equation}
\label{eq:composite}
  \bfsf{P}^\prime = \star^{-1}[T (\star (A\bfsf{  P }\tilde{A})) \tilde{T}].
\end{equation}
\end{rem}

\section{A general transformation} 
\label{sec.gen.trans}

The results of the previous section showed that, 
if we want to study transformations other than 
reflection and rotation, it is not enough to consider
a representation of a vector $\vec{p}$ in terms of a Clifford
algebra based only on $\mathcal{C}\ell_{3,0}$ (represented
by $\mathbf{p}^+$) or $\mathcal{C}\ell_{0,3}$ (represented
by $\mathbf{p}^-$). In fact, let us denote 
$V_3^\pm = \operatorname{span}\{\mathbf{e}_1^\pm,\mathbf{e}_2^\pm,\mathbf{e}_3^\pm\}$ 
and $\mathcal{P}^\pm$ the set of points described by paravectors $\mathsf{P}^\pm$ with 
vector part in $V_3^\pm$. 
The approach based on the Clifford
algebra $\mathcal{C}\ell_{3,3}$ showed that, in the 
case of reflection $N$ and rotation $R$, these transformations 
leave the positive and negative sectors of $\mathcal{C}\ell_{3,3}$ invariant, that is, 
\begin{equation}
N \mathcal{P}^+ \tilde{N} \subset \mathcal{P}^+, \quad N \mathcal{P}^- \tilde{N} \subset 
\mathcal{P}^-, 
\quad R \mathcal{P}^+ \tilde{R} \subset \mathcal{P}^+, \quad 
R \mathcal{P}^- \tilde{R} \subset \mathcal{P}^-  . 
\end{equation}
However, for hyperbolic rotation $H$, shear transformation $S$, 
non-uniform scale transformation $D$ and translation $T$, we have 
\begin{equation}
\Phi \mathcal{P}^+ \tilde{\Phi} \subset \mathcal{P}^+ \oplus \mathcal{P}^- \quad 
\Phi \mathcal{P}^- \tilde{\Phi} \subset \mathcal{P}^+ \oplus \mathcal{P}^- , \quad \Phi \in \{H,S,D,T\} . 
\end{equation}
The same happens for cotranslation:
\begin{equation}
\star[T(\star \mathcal{P}^+)\tilde{T}] \subset \mathcal{P}^+ \oplus \mathcal{P}^- , \quad 
\star[T(\star \mathcal{P}^-)\tilde{T}] \subset \mathcal{P}^+ \oplus \mathcal{P}^-  .
\end{equation}

We want now to study the form of a general transformation of points in the present formalism. The above remarks make it clear that we cannot consider points in $\mathcal{P}^+$
or $\mathcal{P}^-$ separately, but in another subset of $\mathcal{P}^+ \oplus 
\mathcal{P}^-$, and we saw this subset is the one of paravectors with 
vector part of the form $\mathbf{p} = (1/2)(\mathbf{p}^+ + \mathbf{p}^-)$. 
As a consequence, we have to look for transformations such that 
the resulting vector part $\langle\bfsf{P}^\prime\rangle_1$ satisfies 
\begin{equation}
\langle\bfsf{P}^\prime\rangle_1 \cdot \mathbf{e}_i^+ = - 
\langle\bfsf{P}^\prime\rangle_1 \cdot \mathbf{e}_i^-, \quad i=1,2,3. 
\end{equation}
We work with 
$\mathbf{e}_i$ and $\mathbf{e}_i^\ast$ instead of 
$\mathbf{e}_i^+$ and $\mathbf{e}_i^-$, looking for
transformations such that 
\begin{equation}
\label{extra.cond}
\langle\bfsf{P}^\prime\rangle_1 \cdot \mathbf{e}_i = 0 , \quad i = 1,2,3, 
\end{equation}
that is, $\langle\bfsf{P}^\prime\rangle_1$ must not contain any of 
$\mathbf{e}_i^\ast$ ($i=1,2,3$). 
Moreover, we must also take into account the 
condition in eq.\eqref{eq.gar}. 

\begin{thm}
\label{theorem.gen.trans}
Consider the transformation $\bfsf{P} \mapsto \Psi \bfsf{P}\tilde{\Psi}$ 
with $\Psi$ as 
\begin{equation}
\Psi = \Psi_0 + \Psi_1 + \Psi_2 + \Psi_3 + \Psi_4 + \Psi_5 + \Psi_6 , 
\end{equation}
where $\Psi_k$ denotes the $k$-vector part of $\Psi$. Then  
$\Psi \bfsf{P}\tilde{\Psi}$ is a paravector only if 
the following conditions are satisfied: 
\begin{gather}
\label{cond.1}
 2 \Psi_0 \Psi_4 - \Psi_2 \wedge \Psi_2 
-2\Psi_1\wedge\Psi_3 + \Delta_1 = 0 , \\
\label{cond.2}
 2 \Psi_0 \Psi_5 + \Delta_2 = 0 , \\
\label{cond.3}
 2 \Psi_0 \Psi_3 \wedge\mathbf{p} 
- 2\Psi_1\wedge\Psi_2\wedge\mathbf{p} + \Delta_3 = 0 , \\ 
\label{cond.4}
 2 \Psi_0 \Psi_4 \wedge\mathbf{p} 
- \Psi_2\wedge\Psi_2 \wedge\mathbf{p} 
 +2 \Psi_1\wedge\Psi_3\wedge\mathbf{p} 
- 2 \Psi_0 (\mathbf{p}\cdot \Psi_6) + \Delta_4 = 0 , \quad
\end{gather}
where 
\begin{gather}
\label{aux.1}
\Delta_1 =  2\langle \Psi_1 \Psi_5\rangle_4 
+ 2\langle\Psi_2(\Psi_4-\Psi_6)\rangle_4  + 
\langle\Psi_3(-\Psi_3 + 2\Psi_5)\rangle_4
+ \langle \Psi_4 \Psi_4 \rangle_4 , \\
\label{aux.2} 
\Delta_2 = 2\langle \Psi_1(\Psi_4 - \Psi_6)\rangle_5 + 2\langle 
\Psi_2(-\Psi_3+\Psi_5)\rangle_5   + 2 \langle \Psi_3\Psi_4\rangle_5 , \\
\label{aux.3}
\begin{split}
\Delta_3 & =  2\langle\Psi_1\mathbf{p}( \Psi_4 - \Psi_6)\rangle_4 
+ 2\langle\Psi_2\mathbf{p}(-\Psi_3+\Psi_5)\rangle_4  \\
& \quad + 2\langle\Psi_3 \mathbf{p}(\Psi_4-\Psi_6)\rangle_4 + 
2\langle\Psi_4\mathbf{p}\Psi_5\rangle_4 , 
\end{split} \\
\label{aux.4}
\Delta_4  = 
2\langle\Psi_1\mathbf{p}\Psi_5\rangle_5 + 
2 \langle\Psi_2\mathbf{p}(\Psi_4-\Psi_6)\rangle_5 + 
\langle\Psi_3\mathbf{p}(-\Psi_3 + 2\Psi_5)\rangle_5 + 
\langle\Psi_4 \mathbf{p}\Psi_4\rangle_5 .  
\end{gather}
\end{thm}

\medskip

\begin{proof} 
Since $\Psi \bfsf{P}\tilde{\Psi}$ 
is such that $\Psi \bfsf{P}\tilde{\Psi} = (\Psi \bfsf{P}\tilde{\Psi})\widetilde{\;}$, 
the quantity $\Psi \bfsf{P}\tilde{\Psi}$ has scalar, vector, $4$-vector and
$5$-vector parts. For its $4$-vector and $5$-vector parts to be
zero, we must have 
\begin{equation}
\langle \Psi \tilde{\Psi} \rangle_4 = 0, \quad 
\langle \Psi \tilde{\Psi} \rangle_5 = 0 , \quad 
\langle \Psi \mathbf{p} \tilde{\Psi} \rangle_4 = 0 , \quad 
\langle \Psi \mathbf{p} \tilde{\Psi} \rangle_5 = 0 . 
\end{equation}
After some straightforward algebraic manipulations, 
the condition $\langle \Psi \tilde{\Psi} \rangle_4 = 0$ gives 
eq.\eqref{cond.1}, $\langle \Psi \tilde{\Psi} \rangle_5 = 0$ gives
eq.\eqref{cond.2}, $\langle \Psi \mathbf{p} \tilde{\Psi} \rangle_4 = 0$ 
gives eq.\eqref{cond.3} and $\langle \Psi \mathbf{p} \tilde{\Psi} \rangle_5 = 0 $ 
gives eq.\eqref{cond.4}. 
\end{proof}

\medskip

Studying all solutions of this system of equations 
is not a simple task, but, for our interests, studying all solutions is 
not a necessary task. What we are interested in is analysing 
the type of transformations that can be 
described by our formalism. Among the possible transformations, 
let us limit ourselves to the case of continuous transformations, 
and let us suppose $\Psi$ to be close to the identity 
transformation. Our analysis is made easier in this case
because a finite continuous transformation can be obtained by the
repeated application of infinitesimal transformations. Moreover, 
we can consider infinitesimal transformations of the form 
\begin{equation}
\label{infinitesimal.transf}
\Phi_{k} = 1 + \epsilon \, \psi_k , 
\end{equation}
where $\psi_k$ is a $k$-vector and $\epsilon$ is a 
small real parameter. Compositions of 
$\Phi_{k}$ with different $k$-vectors give a infinitesimal
transformation of the form $\Phi = 1 + \epsilon \phi$, 
with $\phi$ a multivector.  

\begin{thm}
\label{theorem.infinitesimal}
Consider the transformation $\bfsf{P} \mapsto \Phi_k \bfsf{P}\tilde{\Phi}_k$ 
with $\Phi_k$ as in eq.\eqref{infinitesimal.transf}. 
Then $\Phi_k \bfsf{P}\tilde{\Phi}_k$ is a paravector if and only 
if $k=0$, $k=1$, and $k=2$ with $\psi_2$ of 
the form $\psi_2 = \mathbf{a}\wedge\mathbf{b}^\ast$. 
\end{thm}

\medskip

\begin{proof} 
See Appendix~\ref{append} for the proof of this theorem.
\end{proof}

\smallskip

Now let us interpret the nature of the general transformations we have. 
To do this interpretation, we will represent the the transformations in matrix form, which will 
allow us a more direct contact with the traditional formalisms, 
facilitating the interpretation. From Theorem~\ref{theorem.infinitesimal} 
we see that the cases to be analysed are $k=0$, $k=1$ and $k=2$, but the case $k=0$ is trivial, as shown in the proof in the Appendix. 
So, let us start with 
the infinitesimal transformations 
\begin{equation}
\Phi_{1} = 1 +\epsilon \, \mathbf{v} 
\end{equation}
and 
\begin{equation}
\Phi_{2} = 1 + \epsilon\, \mathbf{a}\wedge\mathbf{b}^\ast , 
\end{equation}
and consider their composition $\Psi = \Phi_1\Phi_2$. The 
result of $\bfsf{P}^\prime = \Psi \bfsf{P}\tilde{\Psi}$ is 
\begin{equation}
\bfsf{P}^\prime = 1 + \mathbf{p} + \epsilon (2\mathbf{v} + g(\vec{p},\vec{b}) \mathbf{a}) . 
\end{equation}
If we write $\bfsf{P} = p^\prime{}^0 + p^\prime{}^i \mathbf{e}_i$, we can 
represent the result in terms of components as 
\begin{equation}
\begin{pmatrix} p^\prime{}^0 \\ 
p^\prime{}^1 \\ p^\prime{}^2 \\ p^\prime{}^3 \end{pmatrix} = 
\begin{pmatrix} 
1 & 0 & 0 & 0 \\
\epsilon\,2 v^1 & 1 + \epsilon\,a^1 b^1 & \epsilon\, a^1 b^2 & 
\epsilon\, a^1 b^3 \\
\epsilon\,2 v^2 & \epsilon\,a^2 b^1 & 1+ \epsilon\, a^2 b^2 & 
\epsilon\, a^2 b^3 \\
\epsilon\,2 v^3 & \epsilon\,a^3 b^1 & \epsilon\, a^3 b^2 & 
1+ \epsilon\, a^3 b^3 \end{pmatrix} 
\begin{pmatrix} 1 \\ 
p^1 \\ p^2 \\ p^3 \end{pmatrix}  
\end{equation}
This is the matrix representation of an affine transformation. 

On the other hand, we also have transformations of the form 
$\star[\Psi(\star \bfsf{P})\tilde{\Psi}]$. Using 
\begin{gather}
\star[\mathbf{v}(\star 1) + (\star 1)\mathbf{v}] = 0 , \\
\star[\mathbf{v}(\star \mathbf{p}) + (\star \mathbf{p})\mathbf{v}] = g(\vec{v},\vec{p}) , \\
 \star[(\mathbf{a}\wedge\mathbf{b}^\ast)(\star 1) - (\star 1)(\mathbf{a}\wedge\mathbf{b}^\ast)] = g(\vec{a},\vec{b}) , \\
 \star[(\mathbf{a}\wedge\mathbf{b}^\ast)(\star \mathbf{p}) - (\star \mathbf{p})(\mathbf{a}\wedge\mathbf{b}^\ast)] = 
 g(\vec{a},\vec{b})\mathbf{p} - g(\vec{a},\vec{p})\mathbf{b} , 
 \end{gather}
the result of the transformation $\star[\Psi(\star \bfsf{P})\tilde{\Psi}]$ for 
$\Psi = \Phi_1\Phi_2$ can be written in terms of components as 
 \begin{equation} 
\begin{pmatrix} p^\prime{}^0 \\ 
p^\prime{}^1 \\ p^\prime{}^2 \\ p^\prime{}^3 \end{pmatrix}  = 
\begin{pmatrix} 
1  & \epsilon\, v^1 & \epsilon\, v^2 & \epsilon\, v^3 \\
0 &  1 
 -\epsilon\,a^1 b^1 & -\epsilon\, a^2 b^1 & 
-\epsilon\, a^3 b^1 \\
0& -\epsilon\,a^1 b^2 & 1 -\epsilon\, a^2 b^2 & 
-\epsilon\, a^3 b^2 \\
0 & -\epsilon\,a^1 b^3 & -\epsilon\, a^2 b^3 & 
1 -\epsilon\, a^3 b^3 \end{pmatrix} 
\begin{pmatrix} 1 \\ 
p^1 \\ p^2 \\ p^3 \end{pmatrix}  
+ \epsilon \, g(\vec{a},\vec{b}) \begin{pmatrix} 1 \\ 
p^1 \\ p^2 \\ p^3 \end{pmatrix} 
\end{equation}

The composition of both types of transformations (as in \eqref{eq:composite})
can be viewed as 
\begin{equation}
\begin{pmatrix} p^\prime{}^0 \\[1ex] [p^\prime] \end{pmatrix} = 
\begin{pmatrix} a & [v_1] \\[1ex] [v_2] & [A] \end{pmatrix} \begin{pmatrix} 1 \\[1ex] [p]\end{pmatrix} ,
\end{equation}
where $[p]$, $[p^\prime]$ and $[v_2]$ are $1\times 3$ matrices, 
$[v_1]$ is a $3 \times 1$ matrix and $[A]$ is a $3\times 3$ matrix, and 
which can be interpreted as a general 3D projective transformation.

\section{Conclusions}
\label{sec.7}

We have used the Clifford algebra $\mathcal{C}\ell_{3,3}$ 
to describe some geometric transformations  
relevant for computer graphics. The model 
is based on the identification of the vector
space $V_3$ inside the algebra $\mathcal{C}\ell_{3,3}$ 
through the map $\gamma(\vec{v}) = \mathbf{v} = 
\frac{1}{2}(\mathbf{v}^+ + \mathbf{v}^-) $.
This map is the analog of the natural map discussed
in Section 4 of \cite{VM2018}. We have used  
the ideas of~\cite{VM2018} concerning
the use of paravectors $\bfsf{P}$ to represent 
points to study some of their transformations 
of the form $\bfsf{P} \mapsto \Phi \bfsf{P}\tilde{\Phi}$. 
This study showed that reflections and
rotations leave invariant the subspaces that
can be associated with the algebras 
$\mathcal{C}\ell_{3,0}$ and $\mathcal{C}\ell_{0,3}$. 
However, other transformations like shear or non-uniform scale transformations 
do not leave these subspaces invariant. Therefore, 
reflection and rotation can be  
described by $\mathcal{C}\ell_{3,0}$ and $\mathcal{C}\ell_{0,3}$, 
but not transformations like shear or non-uniform scale transformations. 
By using cotranslation, we obtained the single perspective projected
required in computer graphics.  While using both translation and and
cotranslation limits how much the transformation versors can be
combined (with two being required for the computer graphics setting), we note
that with this approach that in addition to the perspective projection
of computer graphics, we obtain all 3D projective transformations.

Compared with other works in the literature, our approach is closely  
related to the approach by Goldman-Mann (GM)~\cite{GS}. This approach 
is based on the Clifford algebra $\mathcal{C}\ell_{4,4}$---that 
is $R(4,4)$ in the notation of \cite{GS}.\footnote{In \cite{GS} the authors have
used Latin indexes $i,j = 0,1,2,3$ but we have changed to 
Greek indexes $\mu,\nu = 0,1,2,3$ since we have already 
used in this work Latin indexes with values $i,j = 1,2,3$.} According to their 
notation, given a basis $\{w_\mu\}$ ($\mu=0,1,2,3$) of the
vector space $W \simeq \mathbb{R}^4$ and a basis 
$\{w_\mu^\ast\}$ ($\mu = 0,1,2,3$) of the dual vector space 
$W^\ast$, we have 
\begin{equation}
w_\mu\cdot w_\nu = 0 , \qquad w_\mu^\ast\cdot w_\nu^\ast = 0 , \qquad 
w_\mu \cdot w_\nu^\ast = \frac{1}{2}\delta_{\mu\nu} , 
\end{equation}
for $\mu,\nu = 0,1,2,3$. 
In \cite{GS} we also have
\begin{equation}
e_\mu = w_\mu + w_\mu^\ast, \qquad \bar{e}_\mu = w_\mu - w_\mu^\ast, 
\end{equation}
for which 
\begin{equation}
e_\mu \cdot e_\nu = \delta_{\mu\nu} , \qquad \bar{e}_\mu \cdot \bar{e}_\nu = 
-\delta_{\mu\nu} , \qquad e_\mu \cdot \bar{e}_\nu = 0 , 
\end{equation}
for $\mu,\nu = 0,1,2,3$. When we compare these relations with 
the ones in Section~\ref{sec.3}, we see that vectors 
and covectors are interpreted in the same way in both approaches, 
with the noticeable difference that the starting point of 
\cite{GS} is a vector space $W \simeq \mathbb{R}^4$, while
our starting point is a vector space $V_3 \simeq \mathbb{R}^3$, 
which leads to the Clifford algebra $\mathcal{C}\ell_{3,3}$ 
in contrast to the Clifford algebra $\mathcal{C}\ell_{4,4}$ of 
\cite{GS}. When we compare the notations, \cite{GS} uses 
$e_i$ and $\bar{e}_i$ while we use $\mathbf{e}_i^+$ and 
$\mathbf{e}_i^-$, respectively, and \cite{GS} uses 
$w_i$ and $w_i^\ast$ while we use $\mathbf{e}_i$ and 
$\mathbf{e}_i^\ast$, respectively.  

Although our model has one dimension less in the starting vector space 
than the GM model, within the Clifford algebra the 
vectors are modelled in the same way, that is, with a positive 
part and a negative part. The fact that our model does not 
use the vectors $e_0$ and $\bar{e}_0$, that is, 
$\mathbf{e}_0^+$ and $\mathbf{e}_0^-$, respectively, in our notation,  
makes the approach to certain problems different in both models. 
In \cite{GS} the translation operation uses $\mathbf{e}_0^+$ and $\mathbf{e}_0^-$, 
but in our approach they are not needed---the difference being that while 
in \cite{GS} the generator of the translation is a bivector, 
in our approach it is a vector. 
However, the situation concerning the perspective projection, 
which also uses $\mathbf{e}_0^+$ and $\mathbf{e}_0^-$ in \cite{GS}, 
is not so simple. To find a description of the 
perspective projection, 
we had to define and use the concept of cotranslation. 
On the other hand, quadric surfaces, 
and their transformation properties, 
can be described in the model \cite{GS},  but 
their description in our model (along with their transformation 
properties) is still under investigation; here again 
the model \cite{GS} uses $\mathbf{e}_0^+$ and $\mathbf{e}_0^-$ in the  
approach to quadric surfaces. 

In a sense we can say that our model and the GM  one 
are based on the same idea and differ in the 
different ways of implementing it. There are, however, 
other models that use Clifford algebras but
are based on different ideas, and therefore on 
different ways of representing geometric objects
and their transformations. One such model 
was provided by Dorst \cite{Dorst2} and is 
based on the study of oriented projective 
transformations of lines. In \cite{Dorst2} 
an interesting comparison is made between 
the models of Dorst, GM and Klawitter 
\cite{Klawitter}. Because of the
relationship between our model and GM one, 
the comparison made in \cite{Dorst2} also 
serves as a comparison between the model presented here 
and the ones of Dorst and of Klawitter, although the criteria
used in~\cite{Dorst2} are more relevant to projective geometry
than they are to computer graphics.
There is also the PGA model developed by Gunn~\cite{Gunn} 
based on a dual projectivized Clifford algebra. 
The relevant Clifford algebra is the one associated
with the degenerate space $\mathbb{R}^{3,0,1}$, 
where the generators satisfy $(E_i)^2 = 1$ ($i=1,2,3$) 
and $(E_0)^2 = 0$. This approach mixes projective, 
ideal and euclidean elements with nice features,
although the transformations in PGA are limited to
rigid body transformations.

Finally, the comparison of our model in terms of 
$\mathcal{C}\ell_{3,3}$ and its version in 
\cite{VM2018} shows that some algebraic manipulations 
are easier to be done with the Clifford approach. 
This statement is based on the fact that the proofs 
of theorems 1 to 6 (especially those of theorems 2 and 3) 
are simpler when done with Clifford algebras compared 
to the proofs of these theorems made in \cite{VM2018}. 
On the other hand, the proof of theorem 7 is not as 
easy in $\mathcal{C}\ell_{3,3}$ than in \cite{VM2018}, 
the main reason being the fact that the Hodge duality 
operator uses the volume element of $V_3$, which does
not have a simple expression in terms of $\mathcal{C}\ell_{3,3}$ 
as in eq.\eqref{omega.V}. Moreover, to justify the 
key point of our model, which is the identification of the vector
space $V_3$ inside the algebra $\mathcal{C}\ell_{3,3}$ 
through the map $\gamma(\vec{v}) = \mathbf{v} = 
\frac{1}{2}(\mathbf{v}^+ + \mathbf{v}^-) $, 
we need to resort to the natural map of \cite{VM2018}
(Section 4), which in our opinion gives a solid 
interpretation for the map $\gamma(\vec{v})$.

\bigskip

\noindent \textbf{\sffamily Acknowledgements:} JV gratefully 
acknowledges the support of a research grant from 
FAPESP---process 2016/21370-9. SM is grateful for the support of the
Natural Sciences and Engineering Research Council of Canada. 
JV is grateful to the University of Waterloo for the hospitality 
during his stay as visiting professor. 

\appendix
\section{Proof of Theorem 6.2 and Some Examples}
\label{append}

{{\bfseries Theorem~6.2} \em
Consider the transformation $\bfsf{P} \mapsto \Phi_k \bfsf{P}\tilde{\Phi}_k$ 
with $\Phi_k$ as in eq.\eqref{infinitesimal.transf}. 
Then $\Phi_k \bfsf{P}\tilde{\Phi}_k$ is a paravector if and only 
if $k=0$, $k=1$, and $k=2$ with $\psi_2$ of 
the form $\psi_2 = \mathbf{a}\wedge\mathbf{b}^\ast$.
}

\medskip

\begin{proof}
 We need to analyse the 
conditions of Theorem~\ref{theorem.gen.trans} and also the 
condition in eq.\eqref{extra.cond}, which implies that
the quantities $\langle \psi_k + \tilde{\psi}_k \rangle_1$, 
$\langle \psi_k \mathbf{p}+\mathbf{p}\tilde{\psi}_k \rangle_1$, 
$\langle \psi_k \tilde{\psi}_k \rangle_1$ and $
\langle \psi_k \mathbf{p}\tilde{\psi}_k \rangle_1$ 
cannot contains terms involving $\mathbf{e}_i^\ast$ ($i=1,2,3$). 
Let us proceed with each case, keeping in mind to neglect 
the second order terms in $\epsilon$.

\smallskip
\noindent \textit{(i)} For $k=0$ we have $\Psi_0 = 1+\epsilon \psi_0$ 
and $\Psi_j = 0$ for $j\neq 0$, and all conditions of Theorem~\ref{theorem.gen.trans} 
are trivially satisfied, as well as the condition in eq.\eqref{extra.cond}. 
The transformation changes the weight of the point, but
does not change the location of the point, because the vector part is
multiplied by the same scalar quantity. \smallskip \\
\noindent \textit{(ii)} For $k=1$ we have $\Psi_0 = 1$, $\Psi_1 = \epsilon \psi_1$ 
and $\Psi_j = 0$ for $j \neq \{0,1\}$.  All conditions of Theorem~\ref{theorem.gen.trans} 
are trivially satisfied, and the condition in eq.\eqref{extra.cond} implies that 
$\psi_1$ cannot contain terms involving 
$\mathbf{e}_i^\ast$ ($i=1,2,3$). Then $\psi_1 = \mathbf{v}$. \smallskip \\
\noindent \textit{(iii)} For $k=2$ we have $\Psi_0 = 1$, $\Psi_2 = \epsilon \psi_2$ 
and $\Psi_j = 0$ for $j \neq \{0,2\}$. The condition given in eq.\eqref{cond.1} gives 
$\psi_2 \wedge\psi_2 = 0$ and the other conditions are trivially satisfied. The 
condition in eq.\eqref{extra.cond} implies that $\psi_2$ cannot  contain 
terms of the form $\mathbf{e}^\ast_i \wedge\mathbf{e}^\ast_j$. 
Moreover, terms of the form $\mathbf{e}_i \wedge\mathbf{e}_j$ produces 
a null transformation, that is, $ (\mathbf{e}_i \wedge\mathbf{e}_j) \mathbf{p} -
\mathbf{p} (\mathbf{e}_i \wedge\mathbf{e}_j) = 0 $ and $(\mathbf{e}_i \wedge\mathbf{e}_j) 
\mathbf{p} (\mathbf{e}_i \wedge\mathbf{e}_j) = 0$, so we will not 
take it into account. 
The condition $\psi_2\wedge\psi_2=0$ is clearly satisfied when 
$\psi_2$ is of the form 
$\psi_2 = \mathbf{a}\wedge\mathbf{b}^\ast$. 
If we write $\psi_2$ as $\mathbf{a}_1\wedge\mathbf{b}_1^\ast + 
\mathbf{a}_2\wedge\mathbf{b}_2^\ast$ the condition $\psi_2 \wedge\psi_2 =0$ 
gives $\mathbf{a}_1 = c_1 \mathbf{a}_2$ or $\mathbf{b}_1^\ast  = 
c_2 \mathbf{b}_2^\ast $ ($c_1$ and $c_2$ constants), which 
also gives $\psi_2$ of the form $\mathbf{a}\wedge\mathbf{b}^\ast$. \smallskip \\
\noindent \textit{(iv)} For $k=3$ we have $\Psi_0 = 1$, $\Psi_3 = \epsilon \psi_3$ 
and $\Psi_j=0$ for $j\neq \{0,3\}$. The condition in eq.\eqref{cond.1} gives 
$\langle \psi_3 \psi_3 \rangle_4  =0$, the condition in eq.\eqref{cond.3} gives
$\psi_3 \wedge\mathbf{p} = 0$, and the other two conditions are trivially satisfied. 
Since $\mathbf{p}$ is arbitrary, we have $\psi_3 = 0$. \smallskip \\
\noindent \textit{(v)} For $k=4$ we have $\Psi_0 = 1$, $\Psi_4 = \epsilon \psi_4$ 
and $\Psi_j =0$ for $j \neq \{0,4\}$. The condition in eq.\eqref{cond.1} gives
(i) $ 2  \psi_4 + \epsilon \langle \psi_4 \psi_4\rangle_4 = 0$, the 
condition in eq.\eqref{cond.4} gives (ii)
$\text{(ii)} \quad 2  \psi_4 \wedge \mathbf{p} + 
\epsilon \langle \psi_4 \mathbf{p}\psi_4\rangle_5 = 0$, 
 and the other
two conditions are trivially satisfied. A direct calculation gives 
\begin{equation*}
\begin{split}
& \langle \psi_4 \mathbf{p}\psi_4\rangle_5 = \langle \psi_4 \mathbf{p}\mathbf{I}^2 
\psi_4 \rangle_5 = - \langle \psi_4 \mathbf{I}\mathbf{p}\psi_4\mathbf{I}\rangle_5 = 
-\langle \psi_2^{\scriptscriptstyle \maltese}\mathbf{p}\psi_2^{\scriptscriptstyle \maltese}
\rangle_5 \\
& = -\psi_2^{\scriptscriptstyle \maltese} \wedge\mathbf{p}\wedge
\psi_2^{\scriptscriptstyle \maltese} = - \psi_2^{\scriptscriptstyle \maltese}\wedge
\psi_2^{\scriptscriptstyle \maltese} \wedge\mathbf{p} = 
-\langle \psi_2^{\scriptscriptstyle \maltese}\psi_2^{\scriptscriptstyle \maltese}\rangle_4 
\wedge\mathbf{p} = -\langle \psi_4\psi_4\rangle_4\wedge\mathbf{p} , 
\end{split}
\end{equation*}
where we used $\mathbf{I}^2 = 1$, $\mathbf{p}\mathbf{I} = -\mathbf{I}\mathbf{p}$ 
and $\psi_2^{\scriptscriptstyle \maltese} = {\scriptstyle \maltese} \psi_4$ 
(see Subsection~\ref{hodge.subsec}). 
Then the second condition can be written as 
$(2 \psi_4 - \epsilon \langle \psi_4 \psi_4\rangle_4)\wedge\mathbf{p} = 0$. 
From the fact that $\mathbf{p}$ is arbitrary, and comparing the 
two conditions, it follows that $\psi_4 = 0$. \smallskip \\
\noindent \textit{(vi)} For $k=5$ we have $\Psi_0 = 1$, $\Psi_5 = \epsilon \psi_5$ 
and $\Psi_j = 0$ for $j\neq \{0,5\}$. The condition in eq.\eqref{cond.2} gives
$2  \psi_5 = 0$ and the other three conditions are trivially satisfied. 
Therefore we have $\psi_5 = 0$. \smallskip \\
\noindent \textit{(vii)} For $k=6$ we have $\Psi_0 = 1$, $\Psi_6  = \epsilon \psi_6$ 
and $\Psi_j = 0$ for $j\neq \{0,6\}$. The condition in eq.\eqref{cond.4} gives
$2(\mathbf{p}\cdot \psi_6) = 0$, and the other three conditions are
trivially satisfied. From the arbitrariness of $\mathbf{p}$, we conclude
that $\psi_6 = 0$. 
\end{proof}

\medskip

The composition of transformations of different types should 
also satisfy the conditions from Theorem~\ref{theorem.gen.trans}, 
although in this case there may be components with different 
$k$-vector parts.
We now consider three types of transformations 
given by the composition of 
the transformations of eq.\eqref{infinitesimal.transf} with 
$k=1$ and $k=2$
and show that conditions (\ref{cond.1}), (\ref{cond.2}), (\ref{cond.3}) and (\ref{cond.4}) of Theorem~\ref{theorem.gen.trans} are met.

\smallskip

\noindent \textit{Example 1.} Let $\Phi_1 = 1 + \epsilon \mathbf{v}$ and
$\Phi_1^\prime = 1 + \eta \mathbf{u}$. The transformation described 
by $\Phi_1 \Phi_1^\prime$ is such that 
\begin{equation*}
 \Psi_0 = 1, \quad 
\Psi_1 = \epsilon \mathbf{v} + \eta\mathbf{u}, \quad 
\Psi_2 = 
\epsilon\eta \mathbf{v}\wedge\mathbf{u} . 
\end{equation*}
 The condition from 
eq.\eqref{cond.2} is trivially satisfied, and 
from eq.\eqref{cond.1} we 
have $\Psi_2 \wedge\Psi_2 = 0$, from eq.\eqref{cond.3} we 
have $\Psi_1\wedge\Psi_2\wedge\mathbf{p} = 0$ 
and from eq.\eqref{cond.4} we have $\Psi_2\wedge\Psi_2\wedge\mathbf{p} = 0$, 
which are satisfied for $\Psi_1$ and $\Psi_2$ given above. 

\smallskip

\noindent \textit{Example 2.} Let $\Phi_1 = 1 + \epsilon \mathbf{v}$ and
$\Phi_2 = 1 + \eta \mathbf{a}\wedge\mathbf{b}^\ast$. The composed transformation 
$\Phi_1 \Phi_2$ is such that 
\begin{equation*}
\Psi_0 = 1, \quad \Psi_1 = \epsilon\mathbf{v} - 
\frac{1}{2}\epsilon\eta g(\vec{v},\vec{b})\mathbf{a} , \quad  
 \Psi_2 = \eta \mathbf{a}\wedge\mathbf{b}^\ast , \quad \Psi_3 = 
\epsilon\eta \mathbf{v}\wedge\mathbf{a}\wedge\mathbf{b}^\ast . 
\end{equation*}
We have, from 
eq.\eqref{cond.1}, 
\begin{equation*}
\text{(i)}\quad \Psi_2\wedge\Psi_2 + 2\Psi_1\wedge\Psi_3 = 0 , 
\end{equation*}
from eq.\eqref{cond.2}, 
\begin{equation*}
\text{(ii)} \quad  \Psi_2\wedge \Psi_3  = 0 , 
\end{equation*}
from eq.\eqref{cond.3}  
\begin{equation*}
\text{(iii)} \quad \Psi_0 \Psi_3\wedge\mathbf{p} - 
\Psi_1\wedge\Psi_2\wedge\mathbf{p} - \langle \Psi_2\mathbf{p}\Psi_3\rangle_4 = 0 , 
\end{equation*}
and from eq.\eqref{cond.4} 
\begin{equation*}
\text{(iv)} \quad \Psi_2\wedge\Psi_2\wedge\mathbf{p} 
-2 \Psi_1\wedge\Psi_3\wedge\mathbf{p} + \langle \Psi_3\mathbf{p}\Psi_3\rangle_5 = 0 . 
\end{equation*}
Using the expressions for $\Psi_j$ ($j=0,1,2,3$) it is immediate that
$\Psi_2 \wedge\Psi_2 = 0$, $\Psi_1\wedge\Psi_3 = 0$, $\Psi_2\wedge\Psi_3 = 0$, 
so that (i) and (ii) are satisfied, and (iv) reduces to $\langle \Psi_3\mathbf{p}\Psi_3\rangle_5 = 0$. It is also easy to see that $\Psi_0\Psi_3 = \Psi_1\wedge\Psi_2$, 
so that (iii) reduces to $\langle \Psi_2\mathbf{p}\Psi_3\rangle_4 = 0$. 
Using the expressions for $\Psi_2$ and $\Psi_3$ and 
\begin{equation*}
\mathbf{p} ( \mathbf{v}\wedge\mathbf{a}\wedge\mathbf{b}^\ast) = \mathbf{p}\wedge
\mathbf{v}\wedge\mathbf{a}\wedge\mathbf{b}^\ast + 
\frac{1}{2}g(\vec{p},\vec{b}) \mathbf{p}\wedge
\mathbf{v}\wedge\mathbf{a} 
\end{equation*}
 we can conclude  
that $\langle \Psi_2\mathbf{p}\Psi_3\rangle_4 = 0$ 
and $\langle \Psi_3\mathbf{p}\Psi_3\rangle_5 = 0$. 

\smallskip

\noindent \textit{Example 3.} Let $\Phi_2 = 1 + \epsilon \mathbf{u}\wedge
\mathbf{v}^\ast$ and $\Phi_2^\prime = 1 + \eta \mathbf{a}\wedge\mathbf{b}^\ast$. 
The composed transformation $\Phi_2 \Phi_2^\prime$ has
\begin{equation*}
\begin{split}
& \Psi_0 = 1 + \frac{1}{4}\epsilon\eta\, g(\vec{u},\vec{b})g(\vec{v},\vec{a}) , \\
& \Psi_2 = \epsilon \mathbf{u}\wedge\mathbf{v}^\ast + \eta \mathbf{a}\wedge
\mathbf{b}^\ast + \frac{1}{2}\epsilon\eta\,[g(\vec{v},\vec{a})\mathbf{u}\wedge
\mathbf{b}^\ast + g(\vec{u},\vec{b})\mathbf{v}^\ast\wedge\mathbf{a}] , \\
&\Psi_4 = \epsilon\eta\, \mathbf{u}\wedge\mathbf{v}^\ast\wedge\mathbf{a}
\wedge\mathbf{b}^\ast .
\end{split}
\end{equation*}
 Eq.\eqref{cond.1} gives 
\begin{equation*}
\text{(i)} \quad 2\Psi_0\Psi_4 - 
\Psi_2 \wedge\Psi_2 + 2\langle \Psi_2\Psi_4\rangle_4 + \langle\Psi_4 \Psi_4
\rangle_4 = 0 .  
\end{equation*}
Eq.\eqref{cond.4} gives 
\begin{equation*}
\text{(ii)} \quad 2\Psi_0\Psi_4\wedge\mathbf{p} - 
\Psi_2 \wedge\Psi_2\wedge\mathbf{p} + 2\langle \Psi_2\mathbf{p}\Psi_4\rangle_5 + 
\langle\Psi_4 \mathbf{p}\Psi_4 \rangle_5 = 0 , 
\end{equation*}
 and the two other conditions
are trivially satisfied. From the expression for $\Psi_0$, $\Psi_2$ and $\Psi_4$ 
we see that $2\Psi_0\Psi_4 - \Psi_2\wedge\Psi_2 = 0$, so conditions (i) and (ii) 
reduce to 
\begin{equation*}
\text{(iii)} \quad 2\langle \Psi_2\Psi_4\rangle_4 + \langle\Psi_4 \Psi_4
\rangle_4 = 0 
\end{equation*}
 and 
\begin{equation*}
\text{(iv)} \quad 2\langle \Psi_2\mathbf{p}\Psi_4\rangle_5 + 
\langle\Psi_4 \mathbf{p}\Psi_4 \rangle_5 = 0 . 
\end{equation*}
 Let us show that
$\langle \Psi_2 \Psi_4 \rangle_4 = 0$. First we calculate $\langle 
(\mathbf{u}\wedge\mathbf{v}^\ast)(\mathbf{u}\wedge\mathbf{v}^\ast
\wedge\mathbf{a}\wedge\mathbf{b}^\ast)\rangle_4$. Using 
eq.\eqref{clifford product} and eq.\eqref{Leibniz} we obtain 
\begin{equation*}
\langle \mathbf{u}\mathbf{v}^\ast (\mathbf{u}\wedge\mathbf{v}^\ast
\wedge\mathbf{a}\wedge\mathbf{b}^\ast)\rangle_4 = 
\frac{1}{2}g(\vec{v},\vec{u}) (\mathbf{u}\wedge\mathbf{v}^\ast
\wedge\mathbf{a}\wedge\mathbf{b}^\ast) 
\end{equation*}
 and 
\begin{equation*}
\langle \mathbf{v}^\ast \mathbf{u}(\mathbf{u}\wedge\mathbf{v}^\ast
\wedge\mathbf{a}\wedge\mathbf{b}^\ast)\rangle_4 = 
\frac{1}{2}g(\vec{v},\vec{u}) (\mathbf{u}\wedge\mathbf{v}^\ast
\wedge\mathbf{a}\wedge\mathbf{b}^\ast) ,
\end{equation*}
and for 
$\mathbf{u}\wedge\mathbf{v}^\ast = (\mathbf{u}\mathbf{v}^\ast
-\mathbf{v}^\ast\mathbf{u})/2$ we obtain 
\begin{equation*}
\langle 
(\mathbf{u}\wedge\mathbf{v}^\ast)(\mathbf{u}\wedge\mathbf{v}^\ast
\wedge\mathbf{a}\wedge\mathbf{b}^\ast)\rangle_4 = 0 . 
\end{equation*} 
Analogously we have 
\begin{equation*}
\begin{split}
& \langle 
(\mathbf{a}\wedge\mathbf{b}^\ast)(\mathbf{u}\wedge\mathbf{v}^\ast
\wedge\mathbf{a}\wedge\mathbf{b}^\ast)\rangle_4 = 0 , \\
& \langle 
(\mathbf{u}\wedge\mathbf{b}^\ast)(\mathbf{u}\wedge\mathbf{v}^\ast
\wedge\mathbf{a}\wedge\mathbf{b}^\ast)\rangle_4 = 0 , \\ 
& \langle 
(\mathbf{a}\wedge\mathbf{v}^\ast)(\mathbf{u}\wedge\mathbf{v}^\ast
\wedge\mathbf{a}\wedge\mathbf{b}^\ast)\rangle_4 = 0 , 
\end{split}
\end{equation*}
and then $\langle \Psi_2\Psi_4 \rangle_4 = 0$. The proof 
that $\langle \Psi_2 \mathbf{p}
\Psi_4\rangle_5 = 0$ is analogous. We have 
\begin{equation*}
\langle \mathbf{u}\mathbf{v}^\ast \mathbf{p}(\mathbf{u}\wedge\mathbf{v}^\ast
\wedge\mathbf{a}\wedge\mathbf{b}^\ast)\rangle_5 = 
\frac{1}{2}g(\vec{v},\vec{u}) (\mathbf{p}\wedge\mathbf{u}\wedge\mathbf{v}^\ast
\wedge\mathbf{a}\wedge\mathbf{b}^\ast)
\end{equation*}
 and 
\begin{equation*}
\langle \mathbf{v}^\ast \mathbf{u}\mathbf{p}(\mathbf{u}\wedge\mathbf{v}^\ast
\wedge\mathbf{a}\wedge\mathbf{b}^\ast)\rangle_5 = 
\frac{1}{2}g(\vec{v},\vec{u}) (\mathbf{p}\wedge\mathbf{u}\wedge\mathbf{v}^\ast
\wedge\mathbf{a}\wedge\mathbf{b}^\ast) , 
\end{equation*}
 and therefore 
\begin{equation*}
\langle 
(\mathbf{u}\wedge\mathbf{v}^\ast)\mathbf{p}(\mathbf{u}\wedge\mathbf{v}^\ast
\wedge\mathbf{a}\wedge\mathbf{b}^\ast)\rangle_5 = 0 . 
\end{equation*}
The same holds for similar terms, and consequently $\langle \Psi_2 \mathbf{p}
\Psi_4\rangle_5 = 0$. For the other two terms, because 
$\langle \Psi_4 \mathbf{p}\Psi_4 \rangle_5 = 
-\langle\Psi_4 \Psi_4\rangle_4\wedge\mathbf{p}$, we just need to
consider $\langle\Psi_4 \Psi_4\rangle_4 $, and we will 
show that $\langle\Psi_4 \Psi_4\rangle_4 = 0$. We can do this using 
\begin{equation*}
\begin{split}
\mathbf{u}\wedge\mathbf{v}^\ast
\wedge\mathbf{a}\wedge\mathbf{b}^\ast  = (1/4)& [\mathbf{u}(\mathbf{v}^\ast
\wedge\mathbf{a}\wedge\mathbf{b}^\ast) - 
\mathbf{v}^\ast(\mathbf{u}\wedge\mathbf{a}\wedge\mathbf{b}^\ast) \\
&  
+ \mathbf{a}(\mathbf{u}\wedge\mathbf{v}^\ast
\wedge\mathbf{b}^\ast) - \mathbf{b}^\ast(\mathbf{u}\wedge\mathbf{v}^\ast 
\wedge\mathbf{a})] 
\end{split}
\end{equation*}
 and calculating four terms like 
$\langle \mathbf{u}(\mathbf{v}^\ast \wedge\mathbf{a}\wedge\mathbf{b}^\ast)\mathbf{u}\wedge\mathbf{v}^\ast
\wedge\mathbf{a}\wedge\mathbf{b}^\ast\rangle_4$, 
for example. To calculate $\langle \mathbf{u}(\mathbf{v}^\ast \wedge\mathbf{a}\wedge\mathbf{b}^\ast)\mathbf{u}\wedge\mathbf{v}^\ast
\wedge\mathbf{a}\wedge\mathbf{b}^\ast\rangle_4$, we proceed as above and obtain 
\begin{equation*}
\begin{split}
& \langle (\mathbf{u}\mathbf{v}^\ast\mathbf{a}\mathbf{b}^\ast)\mathbf{u}\wedge\mathbf{v}^\ast
\wedge\mathbf{a}\wedge\mathbf{b}^\ast\rangle_4 = \frac{1}{4}g(\vec{u},\vec{b})
g(\vec{v},\vec{a})\mathbf{u}\wedge\mathbf{v}^\ast
\wedge\mathbf{a}\wedge\mathbf{b}^\ast , \\
&  \langle (\mathbf{u}\mathbf{a}\mathbf{b}^\ast\mathbf{v}^\ast)\mathbf{u}\wedge\mathbf{v}^\ast
\wedge\mathbf{a}\wedge\mathbf{b}^\ast\rangle_4 = \frac{1}{4}[g(\vec{v},\vec{u})g(\vec{a},\vec{b}) -g(\vec{v},\vec{a})
g(\vec{u},\vec{b})] \mathbf{u}\wedge\mathbf{v}^\ast
\wedge\mathbf{a}\wedge\mathbf{b}^\ast , \\
& \langle (\mathbf{u}\mathbf{b}^\ast\mathbf{v}^\ast\mathbf{a})\mathbf{u}\wedge\mathbf{v}^\ast
\wedge\mathbf{a}\wedge\mathbf{b}^\ast\rangle_4 =
-\frac{1}{4}g(\vec{v},\vec{u})
g(\vec{a},\vec{b})\mathbf{u}\wedge\mathbf{v}^\ast
\wedge\mathbf{a}\wedge\mathbf{b}^\ast , \\ 
& \langle (\mathbf{u}\mathbf{b}^\ast\mathbf{a}\mathbf{v}^\ast)\mathbf{u}\wedge\mathbf{v}^\ast
\wedge\mathbf{a}\wedge\mathbf{b}^\ast\rangle_4 =
\frac{1}{4}g(\vec{v},\vec{u})
g(\vec{a},\vec{b})\mathbf{u}\wedge\mathbf{v}^\ast
\wedge\mathbf{a}\wedge\mathbf{b}^\ast , \\
&  \langle (\mathbf{u}\mathbf{a}\mathbf{v}^\ast\mathbf{b}^\ast)\mathbf{u}\wedge\mathbf{v}^\ast
\wedge\mathbf{a}\wedge\mathbf{b}^\ast\rangle_4 =
\frac{1}{4}[g(\vec{v},\vec{u})g(\vec{a},\vec{b})-g(\vec{v},\vec{a})
g(\vec{u},\vec{b})]\mathbf{u}\wedge\mathbf{v}^\ast
\wedge\mathbf{a}\wedge\mathbf{b}^\ast , \\
& \langle (\mathbf{u}\mathbf{v}^\ast\mathbf{b}^\ast\mathbf{a})\mathbf{u}\wedge\mathbf{v}^\ast
\wedge\mathbf{a}\wedge\mathbf{b}^\ast\rangle_4 =
\frac{1}{4}g(\vec{v},\vec{a})
g(\vec{u},\vec{b})\mathbf{u}\wedge\mathbf{v}^\ast
\wedge\mathbf{a}\wedge\mathbf{b}^\ast ,
\end{split}
\end{equation*}
 and collecting these terms within the 
 expression for $\mathbf{u}(\mathbf{v}^\ast \wedge\mathbf{a}\wedge\mathbf{b}^\ast)$ 
we obtain 
\begin{equation*} 
\langle \mathbf{u}(\mathbf{v}^\ast \wedge\mathbf{a}\wedge\mathbf{b}^\ast)\mathbf{u}\wedge\mathbf{v}^\ast
\wedge\mathbf{a}\wedge\mathbf{b}^\ast\rangle_4 = 0 . 
\end{equation*}
 The same holds for the other similar terms, from 
which we conclude that $\langle \Psi_4 \Psi_4\rangle_4 = 0$.

\end{document}